\documentclass[reqno]{amsart}
\usepackage{aliascnt}
\usepackage{hyperref}
\usepackage{amsfonts}
\usepackage{mathrsfs}
\usepackage{amsmath}
\usepackage{amsmath,graphics}
\usepackage{amssymb}
\usepackage{indentfirst,latexsym,bm,amsmath,pstricks,amssymb,
amsthm,graphicx,fancyhdr,float,color}
\usepackage[all]{xy}
\usepackage{amsthm}
\usepackage{amscd}
\usepackage{hyperref}
\usepackage[all]{hypcap}
\usepackage{multirow}

\hypersetup{backref,
colorlinks=true}

  \renewenvironment{thebibliography}[1]{
    \begin{oldthebibliography}{#1}
      \setlength{\parskip}{0ex}
      \setlength{\itemsep}{0ex}
  }
  {
    \end{oldthebibliography}
  }

\begin{document}

\newtheorem{thm}{Theorem}   
\def\thetheorem{\unskip}
\newtheorem{prop}[thm]{Proposition}
\def\theproposition{\unskip}
\newtheorem{conj}[thm]{Conjecture}
\def\theconjecture{\unskip}
\newtheorem{cor}[thm]{Corollary}
\newtheorem{lem}[thm]{Lemma}
\newtheorem{sublemma}[thm]{Sublemma}
\newtheorem{observation}[thm]{Observation}
\def\thelemma{\unskip}
\theoremstyle{definition}
\newtheorem{defn}{Definition}
\def\thedefinition{\unskip}
\newtheorem{notation}[defn]{Notation}
\newtheorem{rem}[defn]{Remark}
\def\theremark{\unskip}
\newtheorem{ques}[defn]{Question}
\newtheorem{questions}[defn]{Questions}
\def\thequestion{\unskip}
\newtheorem{exmp}[defn]{Example}
\def\theexample{\unskip}
\newtheorem{problem}[defn]{Problem}
\newtheorem{exercise}[defn]{Exercise}
\numberwithin{thm}{section} \numberwithin{defn}{section}

\def\hal{\unskip\nobreak\hfil\penalty50\hskip10pt\hbox{}\nobreak
\hfill\vrule height 5pt width 6pt depth 1pt\par\vskip 2mm}

\renewcommand{\labelenumi}{(\roman{enumi})}
\newcommand{\A}{\mathfrak{A}}
\newcommand{\mfb}{\mathfrak{b}}
\newcommand{\mfC}{\mathfrak{C}}
\newcommand{\mfD}{\mathfrak{D}}
\newcommand{\mfe}{\mathfrak{e}}
\newcommand{\mff}{\mathfrak{f}}
\newcommand{\mfg}{\mathfrak{g}}
\newcommand{\gl}{\mathfrak {gl}}
\newcommand{\mfh}{\mathfrak{h}}
\newcommand{\mfl}{\mathfrak{l}}
\newcommand{\mfm}{\mathfrak{m}}
\newcommand{\mfp}{\mathfrak{p}}
\newcommand{\mfS}{\mathfrak{S}}
\newcommand{\so}{\mathfrak{so}}
\newcommand{\ssp}{\mathfrak{sp}}
\newcommand{\mfsl}{\mathfrak{sl}}
\newcommand{\mft}{\mathfrak{t}}
\newcommand{\mfX}{\mathfrak{X}}
\newcommand{\mfz}{\mathfrak{z}}
\newcommand{\mfu}{\mathfrak{u}}
\newcommand{\mbA}{\mathbb{A}}
\newcommand{\mbE}{\mathbb{E}}
\newcommand{\mFp}{\mathbb{F}_{p}}
\newcommand{\mFq}{\mathbb{F}_{q}}
\newcommand{\mbG}{\mathbb{G}}
\newcommand{\mbP}{\mathbb{P}}
\newcommand{\mbR}{\mathbb{R}}
\newcommand{\mbV}{\mathbb{V}}
\newcommand{\mbZ}{\mathbb{Z}}
\newcommand{\mcA}{\mathcal{A}}
\newcommand{\mcB}{\mathcal{B}}
\newcommand{\mcC}{\mathcal{C}}
\newcommand{\mcE}{\mathcal{E}}
\newcommand{\mcF}{\mathcal{F}}
\newcommand{\mcG}{\mathcal{G}}
\newcommand{\mcH}{\mathcal{H}}
\newcommand{\mcK}{\mathcal{K}}
\newcommand{\mcM}{\mathcal{M}}
\newcommand{\mcO}{\mathcal{O}}
\newcommand{\mcP}{\mathcal{P}}
\newcommand{\mcQ}{\mathcal{Q}}
\newcommand{\mcR}{\mathcal{R}}
\newcommand{\mcS}{\mathcal{S}}
\newcommand{\mcU}{\mathcal{U}}
\newcommand{\mcZ}{\mathcal{Z}}
\newcommand{\msC}{\mathscr{C}}
\newcommand{\msE}{\mathscr{E}}
\newcommand{\msH}{\mathscr{H}}
\newcommand{\mcJ}{\mathscr{J}}
\newcommand{\msN}{\mathscr{N}}
\newcommand{\msR}{\mathscr{R}}
\newcommand{\msS}{\mathscr{S}}
\newcommand{\msW}{\mathscr{W}}
\newcommand{\B}{\mathrm{B}}
\newcommand{\C}{\mathrm{C}}
\newcommand{\cx}{\mathrm{cx}}
\newcommand{\mrcom}{\mathrm{Com}}
\newcommand{\D}{\mathrm{D}}
\newcommand{\mrdim}{\mathrm{dim}}
\newcommand{\srk}{\mathrm{srk}}
\newcommand{\mrlie}{\mathrm{Lie}}
\newcommand{\mrlt}{\mathrm{LT}}
\newcommand{\mmax}{\mathrm{max}}
\newcommand{\mrmax}{\mathrm{Max}}
\newcommand{\mrmin}{\mathrm{min}}
\newcommand{\nil}{\mathrm{nil}}
\newcommand{\mrP}{\mathrm{P}}
\newcommand{\mrpr}{\mathrm{pr}}
\newcommand{\proj}{\mathrm{Proj}}
\newcommand{\rk}{\mathrm{rk}}
\newcommand{\mrrk}{\mathrm{rk_{ss}}}
\newcommand{\mrrkp}{\mathrm{rk_{p}}}
\newcommand{\mrstab}{\mathrm{Stab}}
\newcommand{\mrrad}{\mathrm{rad}}
\newcommand{\mrsp}{\mathrm{Span}}
\newcommand{\mrV}{\mathrm{V}}
\newcommand{\set}[1]{\left\{#1\right\}}
\newcommand{\HHG}{H^{\cdot}(G,\Bk)}
\newcommand{\HHE}{H^{\cdot}(E,\Bk)}
\newcommand{\HHcG}{H^{\cdot}(\mcG,\Bk)}
\newcommand{\Bk}{\Bbbk}
\newcommand{\iotaE}{\iota_{*,\mcE}(V_{r}(\mcE))}
\newcommand{\ck}{\mathrm{char}(\Bk)}
\newcommand{\Ap}{\Bk{\mbZ/p\mbZ}}
\newcommand{\ug}{\underline{\mfg}}
\newcommand{\cnil}{\mfC^{\nil}}
\newcommand{\Phir}{\Phi^{\mrrad}}
\newcommand{\rrm}{r_{\mathrm{max}}}
\newcommand{\rsm}{r_{\mathrm{smax}}}

\def\Ab{\operatorname{Ab}\nolimits}
\def\Ad{\operatorname{Ad}\nolimits}
\def\ad{\operatorname{ad}\nolimits}
\def\d{\operatorname{\bf{d}}\nolimits}
\def\deg{\operatorname{deg}\nolimits}
\def\ev{\operatorname{ev}\nolimits}
\def\exp{\operatorname{exp}\nolimits}
\def\Ens{\operatorname{Ens}\nolimits}
\def\Ext{\operatorname{Ext}\nolimits}
\def\GL{\operatorname{GL}\nolimits}
\def\Gr{\operatorname{Gr}\nolimits}
\def\gr{\operatorname{gr}\nolimits}
\def\id{\operatorname{id}\nolimits}
\def\Im{\operatorname{Im}\nolimits}
\def\ker{\operatorname{ker}\nolimits}
\def\im{\operatorname{im}\nolimits}
\def\Ker{\operatorname{Ker}\nolimits}
\def\Mat{\operatorname{Mat}\nolimits}
\def\mult{\operatorname{mult}\nolimits}
\def\odd{\operatorname{odd}\nolimits}
\def\ord{\operatorname{ord}\nolimits}
\def\PSL{\operatorname{PSL}\nolimits}
\def\Rad{\operatorname{Rad}\nolimits}
\def\res{\operatorname{res}\nolimits}
\def\SL{\operatorname{SL}\nolimits}
\def\SO{\operatorname{SO}\nolimits}
\def\SP{\operatorname{SP}\nolimits}
\def\Spec{\operatorname{Spec}\nolimits}
\def\Specm{\operatorname{Specm}\nolimits}
\def\Sub{\operatorname{Sub}\nolimits}
\def\T{\operatorname{T}\nolimits}
\def\tr{\operatorname{tr}\nolimits}

\newenvironment{changemargin}[1]{%
  \begin{list}{}{%
    \setlength{\topsep}{0pt}%
    \setlength{\topmargin}{#1}%
    \setlength{\listparindent}{\parindent}%
    \setlength{\itemindent}{\parindent}%
    \setlength{\parsep}{\parskip}%
  }%
  \item[]}{\end{list}}

\parindent=0pt
\addtolength{\parskip}{0.5\baselineskip}

\title{Saturation rank for finite group schemes: Finite groups and Infinitesimal group 
schemes}
\author{Yang Pan}
\address{Mathematisches Seminar, Christian-Albrechts-Universit\"{a}t zu Kiel, 24098 Kiel, Germany} \email{ypan@outlook.de}

\pagestyle{plain}
\begin{abstract}
We investigate the saturation rank of a finite group scheme, defined over an algebraically closed
field $\Bk$ of positive characteristic $p$. We begin by exploring the saturation rank for finite groups
and infinitesimal group schemes. Special attention is given to reductive Lie algebras and the second 
Frobenius kernel of the algebraic group $\SL_{n}$.
\end{abstract}
\maketitle

\section{Introduction}
This paper is concerned with the saturation rank of finite group schemes $\mcG$
that are defined over an algebraically closed field $\Bk$ of characteristic $p>0$. 
The Paper \cite{FS} written by Friedlander and Suslin enables us to consider the 
cohomological support variety $V_{\mcG}$, defined as the maximal ideal spectrum of the
even cohomological ring. The classical contexts of $V_{\mcG}$ concern
a recent study on the representation theory of finite groups G and finite dimensional 
restricted Lie algebras $\mfg$. In the field of finite groups, by virtue of Quillen's 
work, it was shown that the dimension of $V_{G}$ is the maximal rank of an elementary 
abelian $p$-subgroup. When it turns to $\mfg$, by setting 
$\underline{\mfg}=\Spec(U_{0}(\mfg)^{*})$
where $U_{0}(\mfg)$ is the restricted enveloping algebra, we are informed that
$V_{\underline{\mfg}}$ and the restricted nullcone $V(\mfg)$ are naturally homeomorphic 
varieties \cite[1.6,5.11]{SFB2}.
Inspired by the approach of elementary abelian $p$-groups to a finite group and the 
generalized 
definition of elementary abelian subgroup schemes given in \cite{Far4}, we now consider
\begin{align*}
r_{p}(\mcG):=\mmax\set{\cx_{\mcE}(\Bk)\; ;\; \mcE \;\mbox{is an elementary abelian 
subgroup of}\;\mcG}
\end{align*}
where $\cx_{\mcE}(\Bk)$ is the complexity of $\mcE$. By concering the irreducible 
components of 
$V_{\mcG}$ in conjunction with the irreducibility of $V_{\mcE}$, we derive the inequality
$r_{p}(\mcG) \leq \dim V_{\mcG}$.
We probably have found that, the behavior of finite groups and infinitesimal group schemes are rather 
different when applying them to the aforementioned formula. Incidentally it is quite clear if we look
at the case of $\mfg=\mfsl_{2}(\Bk)$, indicating that it may be of no hope to investigate
the variety $V_{\mcG}$ by looking at pieces coming from elementary abelian subgroups
when $\mcG$ is an infinitesimal group scheme.

We are now drawn the attention to the saturation rank defined for finite group schemes.
In \cite{Far3}, the saturation rank $\srk(\mcG)$ is exploited for conditioning the 
indecomposibility of Carlson modules. Though it is defined using the theory of
$p$-points, expounded by Friedlander and Pevtsova in their recent paper \cite{FP}, the
homeomorphism between the space $\mrP(\mcG)$ of $p$-points and the projectivization
of $V_{\mcG}$ gives us an interchangeable interpretation. From its definition, it is readily seen that
\begin{align*}
  \srk(\mcG) \leq r_{p}(\mcG) \leq \dim V_{\mcG}.
\end{align*}
The question now was, how the number $\srk(\mcG)$ reveals the properties of $V_{\mcG}$ or $\mcG$.
For instance, when $\srk(\mcG)=\dim V_{\mcG}$ we have $V_{\mcG}$ is equi-dimensional and there are only
finitely many elementary subgroups of $\mcG$ with complexity equaling $\srk(\mcG)$.
The problem was turned into an investigation of the space $\mrP(\mcG)$.
We have already known that, the consideration of $\mrP(\mcG)$ for a finite group scheme generalizes the earlier 
version of the rank variety for a finite group and the variety of 1-parameter subgroups for
an infinitesimal group scheme. 
It is worth detecting the rank $\srk(\mcG)$ via these two precursors 
as our preliminary exploration, even though it is typically difficult to capture those infinitesimal 
group schemes of higher height.

The purpose of this paper is to study the rank $\srk(\mcG)$ for finite groups, finite dimensional 
restricted Lie algebras and a specific example $\SL_{n(2)}$ arising from the second Frobenius kernel of
the algebraic group $\SL_{n}$. The methods we use range from Quillen Stratification for finite groups, 
nilpotent orbits for reductive Lie algebras
and nilpotent commuting variety for $\SL_{n(2)}$. We achieve this by translating from
the homeomorphism 
\begin{align*}
\xymatrix{ \Psi_{\mcG}: \mrP(\mcG) \ar@{->}[r] & \proj V_{\mcG}}.
\end{align*}
When $\mcG$ is a finite group(or equivalently a constant finite group scheme), the proof of 
Theorem \ref{srk-G}, following \cite{Qui1}\cite{Qui2} shows $\srk(\mcG)$ is the minimal rank 
of a maximal elementary abelian $p$-subgroup. Recall that the dimension of $V_{\mcG}$ is the 
maximal rank $r_{p}(\mcG)$ of an elementary abelian $p$-subgroup. It is straightforwardly seen 
that $\srk(\mcG)$ equals the dimension of $V_{\mcG}$ when $V_{\mcG}$ 
is of equi-dimension, and vice versa. 
When $\mcG$ is an infinitesimal group scheme of height $\leq r$, in virtue of the upper-semicontinuity
of a map defined on the support variety $V_{r}(\mcG)$, we find $\srk(\mcG)$ is determined by
the local data. When it applies to reductive Lie algebras (infinitesimal group schemes of height 
$\leq 1$), the regular nilpotent element will be involved in.
We prove in Theorem \ref{ss-rank}, under certain mild restriction on the reductive algebraic group 
$G$, that $\srk(\mfg)$ coincides with the semisimple rank $\mrrk(G)$ of $G$. 
With being curious, we also consider the higher height case, i.e. the second Frobinus kernel 
$\SL_{n(2)}$. 
As being shown in Theorem \ref{height 2}, the saturation rank $\srk(G_{(r)})$, the height of $G_{(r)}$ 
and the semisimple rank $\mrrk(G)$ might  constitute an equality, giving a generalisation to
all $r$-th Frobenius kernel $G_{(r)}$ of the result of Theorem \ref{ss-rank}.

This paper is organised as follows. In Section 2, we introduce the definition of saturation 
rank for all finite group schemes, specializing to finite groups and infinitesimal group schemes.
The characterisation of the saturation rank for reductive Lie algebras, and 
the geomertric realization using nilpotent orbits for an open set is found in Section 3.
Section 4 deals with the saturation rank for the second Frobenius kernel
of $\SL_{n}$, endeavoring to approach a general result.

{\bf Acknowledgement.} 
The results of this paper are part of the author's doctoral thesis, which he was writing
at the University of Kiel. He would like to thank his advisor, Rolf Farnsteiner, for his continuous support. Furthermore, he thanks the members of his working group for proofreading the paper.

\section{Precursors for finite group schemes}
\subsection{}
$\bf{Notations}.$
In this section, we are to introduce and investigate the saturation rank for
all finite group schemes, with an emphasis on constant finite group schemes 
and infinitesimal group schemes. 
In terms of constant finite group schemes, we show that the saturation rank is determined by their 
maximal elementary abelian subgroup schemes; see Theorem \ref{srk-G}. 
In the context of infinitesimal group schemes, we reveal that the saturation rank is 
controlled by the local data; see Theorem \ref{srk-inf}. 
The techniques we use range from Quillen stratification for finite groups (\cite{Qui1},\cite{Qui2}) 
to support varieties for infinitesimal group schemes (\cite{SFB1}, \cite{SFB2}). 
Before we start, we first recall the definition of the saturation rank proposed by 
Farnsteiner (\cite[Sect. 6.4]{Far3}).

Let $\mcG$ be a finite group scheme. For a subgroup $\mcH\subseteq\mcG$, the canonical 
inclusion map $\xymatrix{\iota_{\mcH}: \mcH\ar@{^{(}->}[r]& \mcG}$ induces a continuous yet not 
necessarily injective map $\xymatrix{\iota_{*,\Bk\mcH}: \mrP(\mcH)\ar@{->}[r] &\mrP(\mcG)}$. The 
definition of $p$-points ensures that 
\begin{align*}
\mrP(\mcG)=\bigcup_{\substack{\mcU \subseteq\mcG \\ \text{unipotent abelian}}}
\iota_{*,\Bk\mcU}(\mrP(\mcU)).
\end{align*}
Motivated by this, we consider the set $\mrmax_{au}(\mathcal{G})$ of maximal abelian unipotent 
subgroups of $\mathcal{G}$ as well as the subsets
\begin{align*}
     \mrmax_{au}(\mathcal{G})_{\ell}:= \{ \mcU\in \mrmax_{au}(\mcG)\; ;\;  \cx_{\mcU}(\Bk)\geq \ell\}
\end{align*}
for every $\ell\geq 1$. Setting $\mrP(\mcG)_{\ell}:=\bigcup_{\mcU\in \mrmax_{au}(\mcG)_{\ell}} 
\iota_{*,\Bk\mcU}(\mrP(\mcU))$, the number
\begin{align*}
      \srk(\mcG):=\mathrm{max}\{ \ell\geq 1 \; ; \; \mathrm{P}(\mcG)=\mathrm{P}(\mcG)_{\ell}\}
\end{align*}
is referred to as the $\emph{saturation rank}$ of $\mcG$.

\begin{rem}
In \cite[Section 6.2.1]{Far4}, the author has proved for any abelian unipotent group scheme 
$\mcU\subseteq \mcG$, there exists a unique elementary abelian subgroup scheme 
$\mcE_{\mcU}\subseteq \mcU$ such that 
$\xymatrix{\iota_{*,\Bk\mcE_{\mcU}}: \mrP(\mcE_{\mcU})\ar@{->}[r]& \mrP(\mcU)}$ is a homeomorphism.
We then consider the set $\mrmax_{ea}(\mathcal{G})$ of maximal elementary abelian subgroups of 
$\mcG$ together with the subsets 
\begin{align*}
\mrmax_{ea}(\mathcal{G})_{\ell}:=\set{\mcE \in \mrmax_{ea}(\mathcal{G})\; ; \; \cx_{\mcE}(\Bk)\geq 
\ell } 
\end{align*}
for every $\ell \geq 1$.
Replacing $\mrmax_{au}(\mathcal{G})_{\ell}$ with $\mrmax_{ea}(\mathcal{G})_{\ell}$ in 
$\mrP(\mcG)_{\ell}$ and redefining $\srk(\mcG)$, we find that,
there is no difference on the number $\srk(\mcG)$ between these two settings assured by the 
homeomorphism $\iota_{*,\Bk\mcE_{\mcU}}$.
\end{rem}

\subsection{}
$\bf{Constant \; finite \; group \; schemes.}$ 
Let $G$ be a finite group. Then it defines a constant functor $\mcG_{G}$ which assign to each 
finitely generated connected commutative $\Bk$-algebra the group $G$ itself. This functor is 
represented by $\Bk^{\times |G|}$, indexed by the elements of $G$, with its $\Bk$-linear dual 
$\Bk G$. We call this finite group scheme $\mcG_{G}$ retrieved from $G$ 
a \emph{constant finite group scheme}.

Let $\mcG_{G}$ be a constant finite group scheme with $G=\mcG_{G}(\Bk)$. We denote by $\HHG$ 
the cohomology ring $H^{*}(G,\Bk)$ of $G$ if $\ck=2$ and the subring $H^{ev}(G,\Bk)$ of elements of 
even degree if $\ck>2$. Evens and Venkov have proved independently that $\HHG$ is a finitely 
generated commutative $\Bk$-algebra. We denote by $V_{G}=\mmax(\HHG)$ the maximal ideal spectrum, 
an affine variety corresponding to $\HHG$. Let $E\leq G$ be an elementary abelian $p$-subgroup of 
$G$. Then there is a restriction map 
\begin{align*}
 \xymatrix{\res_{G,E}: \HHG \ar@{->}[r]  &  \HHE }
\end{align*}
which gives rise to a map of affine varieties $\xymatrix{\res_{G,E}^{*}: V_{E}\ar@{->}[r] & V_{G}}$.
Quillen has investigated the variety $V_{G}$ and shown that it is stratified by pieces coming from 
elementary abelian subgroups of $G$, which is known as Quillen Stratification; see \cite{Qui1} and 
\cite{Qui2} for details. A weak version of his result is the following 
\begin{align*}
V_{G}=\bigcup_{\substack{ E\leq G \\  \text{elemab}}} \res_{G,E}^{*}V_{E},
\end{align*}
where $\res_{G,E}^{*} V_{E}=V(\ker(\res_{G,E}))$ is an irreducible closed subvariety of $V_{G}$.

\begin{notation}
Keep the notation for $\HHG$ and $V_{G}$. We denote by $\HHG^{\dagger}$ the augmentation 
ideal of $\HHG$. Let $I$ be an ideal of $\HHG$, then $\gr(I)$ is defined to be the unique maximal 
homogeneous ideal inside of $I$. We denote by $\proj V_{G}$ the set of homogeneous ideals of 
$\HHG$ which are maximal among
those homogeneous ideals other than the augmentation ideal $\HHG^{\dagger}$. Then $\proj V_{G}$
can be identified with the set of $\gr(\mfm)$ for $\mfm\in V_{G}\setminus \set{\HHG^{\dagger}}$.
Let $E$ be an elementary abelian subgroup of $G$. Observe that
$\res_{G,E}^{*}(\gr(\mfm))=\gr(\res_{G,E}^{*}(\mfm))$ for any $\mfm\in V_{E}\setminus \set{\HHE^{\dagger}}$. 
Then the map $\xymatrix{\res_{G,E}^{*}: V_{E} \ar@{->}[r] & V_{G}}$ induces a map
\begin{align*}
 \xymatrix{\res_{G,E}^{*}: \proj V_{E} \ar@{->}[r] & \proj V_{G}}.
\end{align*}
\end{notation}

An elementary abelian subgroup of $\mcG_{G}$ is isomorphic to $\mcG_{E}$ where $E$ is an elementary 
abelian $p$-subgroup of $G$; see \cite[6.2]{Far4}. Without any real ambiguity, we will use $\mcG_{E}$ and $E$ alternatively 
for the sake of convenience. By denoting
\begin{align*}
V_{G}(\ell)=\bigcup_{\mcG_{E}\in \mrmax_{ea}(\mcG_{G})_{\ell}} \res_{G,E}^{*} V_{E},
\end{align*}
we have the following Lemma:

\begin{lem}{\label{trans1}}
Suppose $\mcG_{G}$ is a constant finite group scheme with $\mcG_{G}(\Bk)=G$. Then
\begin{align*}
\srk(\mcG_{G})=\mmax\set{\ell \geq 1 \; ; \; V_{G}=V_{G}(\ell)}.
\end{align*}
\end{lem}
\begin{proof}
We utilize the homeomorphism $\xymatrix{\Psi_{\mcG_{G}}: \mrP(\mcG_{G})\ar@{->}[r]  & 
\proj V_{G}}$ presented in \cite[Sect. 4]{FP} to verify this. 
First assume that $V_{G}=V_{G}(\ell)$. Let $[\alpha]\in \mrP(\mcG_{G})$ be an equivalence class,
then $\Psi_{\mcG_{G}}([\alpha])\in \proj V_{G}$. 
Thus there is $\mfm\in V_{G}\setminus\set{\HHG^{\dagger}}$ such that 
$\Psi_{\mcG_{G}}([\alpha])=\gr(\mfm)$. By our assumption, there exists
a maximal elementary abelian subgroup scheme $\mcG_{E}$ of $\mcG_{G}$ with 
$\cx_{\mcG_{E}}(\Bk)\geq \ell$ such that $\mfm\in \res_{G,E}^{*}(V_{E}\setminus 
\set{\HHE^{\dagger}})$. Then $\gr(\mfm)\in \res_{G,E}^{*}\proj V_{E}$. 
The bijective map $\Psi_{\mcG_{E}}$ ensures 
$\gr(\mfm)=\res_{G,E}^{*}(\Psi_{\mcG_{E}}([\beta]))=\Psi_{\mcG_{G}}([\iota_{*,\Bk E} \circ \beta])$ for some $[\beta]\in \mrP(\mcG_{E})$.  
As a result, $\Psi_{\mcG_{G}}([\alpha])=\Psi_{\mcG_{G}}([\iota_{*,\Bk E} \circ \beta])$, which 
gives $[\alpha]=\iota_{*, \Bk E}([\beta])$ since $\Psi_{\mcG_{G}}$ is bijective. 
Hence, we have $\mrP(\mcG_{G})=\mrP(\mcG_{G})_{\ell}$.

On the other hand, we assume $\mrP(\mcG_{G})=\mrP(\mcG_{G})_{\ell}$. 
Let $\mfm\in V_{G}$. If $\mfm=\HHG^{\dagger}$, then $\mfm\in \res_{G,E}^{*} V_{E}$ for any maximal 
elementary abelian subgroup $\mcG_{E}$ of $\mcG_{G}$. So it suffices to consider 
$\mfm \neq \HHG^{\dagger} $. Then $\gr(\mfm)=\ker \alpha^{\cdot} $ for some $[\alpha] \in
\mrP(\mcG_{G})$. By our assumption, there is a maximal elementary abelian subgroup $\mcG_{E}$
of $\mcG_{G}$ with $\cx_{\mcG_{E}}(\Bk)\geq \ell$ such that $[\alpha]=\iota_{*,\Bk E}([\beta])$ where 
$[\beta]\in \mrP(\mcG_{E})$. This gives $\ker \alpha^{\cdot}=\ker (\beta^{\cdot}\circ \res_{G,E})$,
and consequently $\ker(\res_{G,E})\subset \gr(\mfm)\subset \mfm$. 
As a result, $\mfm \in V(\ker(\res_{G,E}))=\res_{G,E}^{*} V_{E}$. 
Therefore, we have $V_{G}=V_{G}(\ell)$.

Combining these two, it has been shown that $\mrP(\mcG_{G})=\mrP(\mcG_{G})_{\ell}$ if and only if
$V_{G}=V_{G}(\ell)$ for $\ell \geq 1$. This ultimately gives 
$\srk(\mcG_{G})=\mmax\set{\ell \geq 1 \; ; \; V_{G}=V_{G}(\ell)}$.  
\end{proof}

With the finiteness property of the number of elementary abelian subgroups of a finite group 
behind us, Lemma \ref{trans1} tells us that the saturation rank $\srk(\mcG_{G})$ 
of a constant finite group scheme is clear, i.e. is
the minimal dimension of an irreducible component of $V_{G}$. 
Recall that a maximal elementary abelian subgroup $E$ of $G$ has the property:
\begin{itemize}
\item $E$ is not conjugate to a proper subgroup of any other elementary abelian subgroups.
\end{itemize}
Let $\mcM(G)$ be the set of representatives from each conjugacy class of a maximal elementary
abelian subgroup. The following theorem is dedicated to establishing a relation between the
set $\mcM(G)$ and the set of irreducible components of $V_{G}$.

\begin{thm}{\label{srk-G}}
Let $\mcG_{G}$ be a constant finite group scheme with $\mcG_{G}(\Bk)=G$.
Then the assignment 
\begin{equation*}
 E \longmapsto \res_{G,E}^{*}V_{E}
\end{equation*} 
induces a bijection
\[ 
\left\{
  \begin{gathered} \text{elementary abelian subgroups }\\
\text{$E$ in the set $\mcM(G)$}
\end{gathered}\, \right\} \;\stackrel{\sim}\longrightarrow\; \left\{
\begin{gathered} \text{irreducible components of the }\\ \text{affine variety $V_{G}$}
\end{gathered}\,  \right \}\, .
\] 
\end{thm}
\begin{proof}
We first show that the assignment is well-defined. 
Let $E$ be a maximal elementary abelian subgroup.  We adopt the notions from \cite[Sect. 5.6]{Ben}
as follows:
\begin{align*}
V_{E}^{+}:=V_{E}-\bigcup_{\substack{F<E \\ \text{elemab}}} \res_{E,F}^{*}V_{F}, \qquad
V_{G,E}^{+}:=\res_{G,E}^{*}V_{E}^{+}.
\end{align*}
Then by \cite[Lemma 5.6.2]{Ben}, there exists an element $\varrho_{E}$ of $\HHG$ with the
property $V_{G,E}^{+}=\res_{G,E}^{*}V_{E}-V(\varrho_{E})$.
Suppose additionally that $E^{'}$ is an elementary abelian subgroup which is not conjugated to
$E$. Since $E$ is maximal, we have $\res_{G,E^{'}}(\varrho_{E})=0$ according to 
\cite[Lemma 5.6.2]{Ben}. This subsequently implies $\res_{G,E^{'}}^{*}V_{E^{'}} \subseteq V(\varrho_{E})$, ensuring $\res_{G,E}^{*}V_{E}\nsubseteq \res_{G,E^{'}} V_{E^{'}}$.
Thus, $\res_{G,E}^{*}V_{E}$ is maximal in $V_{G}$, i.e. is an irreducible component.
The bijection follows immediately.
\end{proof}

\begin{cor}
Let $\mcG_{G}$ be a constant finite group scheme with $\mcG_{G}(\Bk)=G$.
Then $\srk(\mcG_{G})$ is the minimal rank of a maximal elementary abelian subgroup.
\end{cor}
\begin{proof}
Since $\dim \res_{G,E}^{*}V_{E}=\dim \HHE=\rk(E)$, the result is readily seen.
\end{proof}

\begin{cor}
Let $\mcG_{G}$ be a constant finite group scheme with $\mcG_{G}(\Bk)=G$. Suppose additionally
$V_{G}$ is equidimensional, then $\srk(\mcG_{G})=\dim V_{G}$, and vice versa.
\end{cor}

\begin{exmp}
We consider the dihedral group $D_{8}$ and suppose $p=2$. It has two generators $a$ and $b$ with  relations:
\begin{align*}
     a^{4}=1=b^{2}, \;\mbox{and}\; \;aba=b.
\end{align*}
There are two maximal elementary abelian $2$-subgroups:
\begin{align*}
  \set{e, a^{2}, b, a^{2}b}, \;\mbox{and}\;\;\set{e, a^{2}, ab, a^{3}b}
\end{align*}
which are isomorphic to $\mbZ/2\mbZ \times \mbZ/2\mbZ$. As a result, 
$\srk(\mcG_{D_{8}})=\dim V_{D_{8}}=2$.
\end{exmp}
\subsection{} 
$\bf{Infinitesimal \; group \; schemes.}$ 
We say a finite group scheme $\mcG$ is \emph{infinitesimal} if its coordinate algebra 
$\Bk[\mcG]$ is a local algebra.
Then the augmentation ideal $\Bk[\mcG]^{\dagger}$ of $\Bk[\mcG]$ is its unique maximal ideal. 
Associated to $\mcG$, it is of height $\leq r\in \mathbb{N}_{0}$ if $x^{p^{r}}=0$ for all 
$x\in \Bk[\mcG]^{\dagger}$. 
A class of infinitesimal group schemes that have served as prototypical examples arise 
from reduced algebraic group schemes $G$ by taking
their $r$-th Frobenius kernel $G_{(r)}$ via the Frobenius map $F^{r}$. 
The representing algebra of $G_{(r)}$ is
then $\Bk[G_{(r)}]=\Bk[G]/I$, where $I$ is an ideal generated by elements $x^{p^{r}}$ for 
$x\in\Bk[G_{(r)}]^{\dagger}$.
In particular, when $r=1$, we find that $\Bk[G_{(1)}]$ is the dual of the restricted enveloping algebra 
of algebraic Lie algebra $\mfg=\mrlie(G)$, a sepecial case of restricted Lie algebras which carries a 
canonical $[p]$-structure equivariant under the adjoint action of $G(\Bk)$. In general, there is a 
categorical equivalence between the category of finite dimensional $p$-restricted Lie algebras and 
category of infinitesimal group schemes of height $\leq 1$. 
Henceforth, for any given finite dimensional restricted Lie algebra $(\mfg,[p])$, we denote the 
associated infinitesimal group scheme $\mcG_{\mfg}:=Spec((U_{0}(\mfg))^{*})$ by $\underline{\mfg}$.

Let $\mcG$ be an infinitesimal group scheme of height $\leq r$. An infinitesimal 1-parameter 
subgroup of $\mcG_{R}$ over a commutative $\Bk$-algebra $R$ is a homomorphism of $R$-group 
schemes $\xymatrix{\mbG_{a(r),R}\ar@{->}[r]&\mcG_{R}}$. Let $V_{r}(\mcG)$ be the functor, which 
sends every commutative $\Bk$-algebra $A$ to the group $V_{r}(\mcG)(A)=Hom_{Gr/A}
(\mbG_{a(r),A},\mcG_{A})$. The functor $V_{r}(\mcG)$ is represented by an affine scheme of finite 
type over $\Bk$, \cite[Theorem 1.5]{SFB1}. The coordinate algebra $\Bk[V_{r}(\mcG)]$ of $V_{r}(\mcG)$ 
is then a graded connected algebra generated by homogeneous elements of degree $p^{i},0\leq i\leq 
r-1$. In what follows, we will only concentrate on the $\Bk$-rational points of the scheme 
$V_{r}(\mcG)$, which is still denoted by $V_{r}(G)$. 
Note that an infinitesimal 
1-parameter subgroup $\xymatrix{v:\mbG_{a(r)}\ar@{->}[r] & \mcG}$ over $\Bk$ may be 
factored as
\begin{displaymath}
\xymatrix{ \mbG_{a(r)} \ar@{->}^{u}[r]  & \mbG_{a(s)} \ar@{^{(}->}[r]  &\mcG } 
\end{displaymath}                         
for some $1\leq s\leq r$, where $\mbG_{a(s)}$ is an 
elementary abelian subgroup (see \cite[6.2]{Far4}).

Let $\mcC_{\mcG}$ be the category whose objects are 
elementary abelian subgroups of $\mcG$, and whose morphisms are inclusions. Similarly, define 
$\mcC_{\Bk\mcG}$ to be the category having commutative Hopf subalgebras of $\Bk\mcG$ whose 
underlying associative algbera is isomorphic to some truncated 
polynomial ring $\Bk[x_{1},\ldots,x_{n}]/(x_{1}^{p},\ldots,x_{n}^{p})$ as its objects, and morphisms are 
also given by inclusions. There is a categorical equivalence between $\mcC_{\mcG}$ and $\mcC_{\Bk
\mcG}$ via the functor $\xymatrix{\mcF:\mcC_{\mcG}\ar@{->}[r] &\mcC_{\Bk\mcG}}$ by sending 
$\mcE$ to $k\mcE$.

Let $\xymatrix{\Gr_{d}(\Bk\mcG): \mrcom_{\Bk} \ar@{->}[r] & \Ens}$ be the Grassmann scheme, i.e., the 
$\Bk$-functor that assign to every commutative $\Bk$-algebra $R$ the set $\Gr_{d}(\Bk\mcG)(R)$ of 
$R$-direct summands of $\Bk\mcG\otimes_{\Bk}R$ of rank $d$; see \cite[Sect I.1.9]{Jan1}. 
We begin with the consideration of the subfunctor $\Sub_{d}(\Bk\mcG)\subseteq \Gr_{d}(\Bk\mcG)$ 
which is given by
\begin{align*}
\Sub_{d}(\Bk\mcG)(R):=\set{V\in \Gr_{d}(\Bk\mcG)(R)\mid V \cdot V\subseteq V, \Delta(V)\subseteq 
    V\otimes_{R} V}
\end{align*}
for every commutative $\Bk$-algebra $R$. Recall that the base change $\Bk\mcG\otimes_{\Bk} R$
of the Hopf $\Bk$-algebra $\Bk\mcG$ is then a Hopf $R$-algebra.

\begin{prop}{\label{Sub}}
Keep the notations for $\Gr_{d}(\Bk\mcG), \Sub_{d}(\Bk\mcG)$ as above. Then the functor
$\Sub_{d}(\Bk\mcG)$ is a closed subfunctor of $\Gr_{d}(\Bk\mcG)$.
\end{prop}
\begin{proof}
Let $\xymatrix{\psi: R\ar@{->}[r] & S}$ be a $\Bk$-algebra homomorphism. Then it induces a Hopf 
$S$-algebra homomorphism $\xymatrix{\id\otimes(\psi\hat{\otimes}\id): (\Bk\mcG\otimes_{\Bk} R)
\otimes_{R}S \ar@{->}[r] & \Bk\mcG\otimes_{\Bk} S}$ by sending $x\otimes r\otimes s$ to $x\otimes 
\psi(r)s$, it follows that 
$\Gr_{d}(\Bk\mcG)(\psi)$ sends $\Sub_{d}(\Bk\mcG)(R)$ to $\Sub_{d}(\Bk\mcG)(S)$. As as result, 
$\Sub_{d}(\Bk\mcG)$ is a subfunctor of the $\Bk$-functor $\Gr_{d}(\Bk\mcG)$.

The closedness of $\Sub_{d}(\Bk\mcG)$ is verified in this fashion: for every commutative $\Bk$-algebra
$A$ and every morphism $\xymatrix{f: \Spec_{\Bk}(A) \ar@{->}[r] & \Gr_{d}(\Bk\mcG)}$,
$f^{-1}(\Sub_{d}(\Bk\mcG))$ is a closed subfunctor of $\Spec_{\Bk}(A)$; see \cite[I.1.12]{Jan1}. 
By Yoneda's Lemma, the morphism corresponds to an $A$-point $W\in \Gr_{d}(\Bk\mcG)(A)$.

Fix a basis $\set{v_{1},\ldots,v_{n}}$ of $\Bk\mcG$. Let $\alpha_{ij\ell}$ be the elements of $A$ that 
are given by
\begin{align*}
  v_{i}\cdot v_{j}=\sum_{\ell=1}^{n} \alpha_{ij\ell}v_{\ell}, \qquad 1\leq i,j\leq n
\end{align*}
Denote by $\set{w_{1},\ldots,w_{d}}$ a set of generators of the locally free $A$-module $W$ of rank 
$d$, and define elements $a_{ri},b_{ri},c_{ri}\in A$
via
\begin{align*}
  & w_{r}=\sum_{i=1}^{n} v_{i}\otimes a_{ri},    1\leq r \leq d,\\
  & \Delta(w_{r})=\sum_{i=1}^{n}\sum_{j=1}^{n}(v_{i}\otimes b_{ri})\otimes (v_{j}\otimes c_{rj}), 1\leq r\leq 
            d.
\end{align*}
By definition of $\Gr_{d}(\Bk\mcG)(A)$, there exists an $A$-submodule $W^{'}\subseteq \Bk\mcG
\otimes_{\Bk}A$
such that 
\begin{align*}
    \Bk\mcG\otimes_{\Bk}A = W\oplus W^{'}
\end{align*}
and we denote by $\xymatrix{\mrpr: \Bk\mcG\otimes_{\Bk}A \ar@{->}[r] & W^{'} }$ the corresponding 
projection. This $A$-linear map is given by 
\begin{align*}
  \mrpr(v_{j}\otimes 1)=\sum_{i=1}^{n}v_{i}\otimes \kappa_{ij},\qquad 1\leq j\leq n.
\end{align*}
We let $I\subseteq A$ be the ideal generated by the elements for $1\leq i,j,t\leq n,1\leq r,s\leq d$
\begin{align*}
h_{rst}=\sum_{\ell=1}^{n}\sum_{i=1}^{n}\sum_{j=1}^{n}\alpha_{ij\ell}a_{ri}a_{sj}\kappa_{t\ell}\;\;\; ;\;
g_{trj}=\sum_{q=1}^{n} b_{rq}\kappa_{tq}c_{rj}, \;\;\; ;\;
\gamma_{tri}=\sum_{q=1}^{n}b_{ri}\kappa_{tq}c_{rq}.
\end{align*}
Now let $\xymatrix{\psi: A \ar@{->}[r] & R}$ be a homomorphism of $\Bk$-algebras. Then we have
\begin{align*}
  \Bk\mcG\otimes_{\Bk} R= W_{R}\oplus W^{'}_{R}
\end{align*}
where $X_{R}=X\otimes_{A} R$ for $X\in\set{W,W^{'}}$. The corresponding projection
\begin{align*} 
\xymatrix{\mrpr_{R}: \Bk\mcG\otimes _{\Bk}R \ar@{->}[r] & W^{'}_{R}}
\end{align*}
is given by
\begin{align*}
  \mrpr_{R}(v_{j}\otimes 1)=\sum_{i=1}^{n} v_{i}\otimes \psi(\kappa_{ij}),\qquad 1\leq j\leq n.
\end{align*}
In view of 
\begin{align*}
  (w_{r}\otimes 1)\cdot (w_{s}\otimes 1)=\sum_{\ell=1}^{n}v_{\ell}\otimes\sum_{i=1}^{n}\sum_{j=1}^{n}
     \alpha_{ij\ell}\psi(a_{ri})\psi(a_{sj})
\end{align*}
this gives rise to
\begin{align*}
 \mrpr_{R}((w_{r}\otimes 1)\cdot (w_{s}\otimes 1))=\sum_{t=1}^{n} v_{t}\otimes \psi(h_{rst}).
\end{align*}
By the same token, we have
\begin{align*}
 (\mrpr_{R}\otimes \id_{R})\circ \Delta_{R}(w_{r}\otimes 1)=\sum_{t=1}^{n}\sum_{j=1}^{n}
 (v_{t}\otimes 1)\otimes (v_{j}\otimes \psi(g_{trj})), \\
 (\id_{R}\otimes \mrpr_{R} )\circ \Delta_{R}(w_{r}\otimes 1)=\sum_{t=1}^{n}\sum_{i=1}^{n}
 (v_{i}\otimes \psi(\gamma_{tri}))\otimes (v_{t}\otimes 1 ).
\end{align*}
Now suppose that $\psi(I)=0$. Then the three forgoing identities imply $W_{R}\cdot W_{R}
\subseteq \ker 
\mrpr_{R}=W_{R}$, as well as $\Delta_{R}(W_{R})\subseteq \ker(\mrpr_{R}\otimes \id_{R})\cap 
\ker(\id_{R}\otimes \mrpr_{R})=(W_{R}\otimes_{R}(\Bk\mcG\otimes_{\Bk}R))\cap((\Bk\mcG\otimes_{\Bk}
R)\otimes_{R}
W_{R})=W_{R}\otimes _{R}W_{R}$, thus $W_{R}\in \Sub_{d}(\Bk\mcG)(R)$. Conversely, if we have 
$W_{R}\in \Sub_{d}(\Bk\mcG)(R)$, then $I \subseteq \ker(\psi)$.
Observe that 
$f^{-1}(\Sub_{d}(\Bk\mcG))(R)=\set{\psi\in\Spec_{\Bk}(A)(R) \mid W_{R}=\Gr_{d}(\Bk\mcG)(\psi)(W)\in 
\Sub_{d}(\Bk\mcG)(R)}$. This shows
\begin{align*}
f^{-1}(\Sub_{d}(\Bk\mcG))(R)=\set{\psi\in\Spec_{\Bk}(A)(R) \mid \psi(I)=0}. 
\end{align*}
Consequently, 
$f^{-1}(\Sub_{d}(\Bk\mcG))$ is a closed subfunctor of $\Spec_{\Bk}(A)$, as desired.
\end{proof}

\begin{rem}{\label{Ab}}
Let $\Ab_{d}(\Bk\mcG)\subseteq \Gr_{d}(\Bk\mcG)$ be the subfunctor of commutative subalgebras
(contains identity element of $\Bk \mcG$) of $\Bk\mcG$. 
It is a closed subfunctor, which can be proved similar to Proposition \ref{Sub}. 
The above proposition shows that the sets $\Sub_{d}(\Bk\mcG)(\Bk), \Ab_{d}(\Bk\mcG)(\Bk)$ of rational 
$\Bk$-points of these functors are closed subsets of the Grassmann variety $\Gr_{d}(\Bk\mcG)(\Bk)$.
\end{rem}

\begin{prop}{\label{proj}}
Let $\mfC(\ell,\mcG)$ be the set consisting of objects of $\mcC_{\mcG}$ having complexity $\ell$. Then $\mfC(\ell,\mcG)$ is a projective variety.
\end{prop}
\begin{proof}
Let $\mathcal{X}:=\set{ H\in \Gr_{p^{\ell}}(\Bk\mcG)(\Bk)\mid H\;\mbox{is a commutative Hopf subalgebra of}
\; \Bk\mcG}$. Then by Lemma 1.2(1) of \cite{Far2} on Hopf algebras, we have
\begin{align*}
\mathcal{X}=\Sub_{p^{\ell}}(\Bk\mcG)(\Bk)\cap \Ab_{p^{\ell}}(\Bk\mcG)(\Bk).
\end{align*}
According to our Remark \ref{Ab}, $\mathcal{X}$ is a closed subvariety of $\Gr_{p^{\ell}}(\Bk\mcG)(\Bk)$. 
Consider the set $\mathcal{Y}:=\mfC(p^{\ell},\Bk\mcG)$, consisting of objects of $\mcC_{\Bk\mcG}$ with dimension 
$p^{\ell}$. Then $\mathcal{Y}\subset \mathcal{X}$ is a subset of $\mathcal{X}$.
By the equivalence of categories between $\mcC_{\mcG}$ and $\mcC_{k\mcG}$, it suffices to endow 
$\mathcal{Y}$ with a projective variety structure, i.e. to show $\mathcal{Y}$ is closed. 

Let $\Bk\mcG_{p}:=\{u\in \Bk\mcG \mid u^{p}=0\}$ be the set of $p$-nilpotent elements in $\Bk\mcG$, then it 
is a closed conical subvariety. By setting $H_{p}:=H\cap \Bk\mcG_{p}$ for $H\in \mathcal{X}$,
we are going to verify that: The underlying associative algebra of $H$ is isomorphic to 
$\Bk[x_{1},\ldots,x_{\ell}]/(x_{1}^{p},\ldots,x_{\ell}^{p})$ if and 
only if $\mathrm{dim}_{\Bk}H_{p} \geq p^{\ell}-1$.
First if we have such algebraic isomorphism for $H$, then $H_{p}=\Rad H$ and $\mathrm{dim}_{\Bk}H_{p}\geq 
p^{\ell}-1$. Conversely, if $\mathrm{dim}_{\Bk}H_{p}\geq p^{\ell}-1$ with $\mathrm{dim}_{\Bk}H=p^{\ell}$, then 
$H=H_{p}\oplus \Bk{1_{H}}$ and $H$ must be local since the identity 
element is the unique non-zero idempotent element in $H$. Notice that $H$ is commutative,
then it represents an infinitesimal group scheme. Theorem in \cite[Sect 14.4]{Wat}
ensures that $H$ has to be isomorphic to a truncated polynomial ring, say 
$\Bk[x_{1},\ldots,x_{t}]/(x_{1}^{p^{e_{1}}},\ldots,x_{t}^{p^{e_{t}}})$.
By the definition of $H$, the $p$-th power of $x_{i}$ for $1\leq i\leq t$ should be zero, i.e. $x_{i}^{p}=0$. This
implies all $e_{i}=1$ and further $t=\sum_{i=1}^{t} e_{i}=\ell$ by dimension. Thus for such $H$, its underlying 
associative algebra is isomorphic to $\Bk[x_{1},\ldots,x_{\ell}]/(x_{1}^{p},\ldots,x_{\ell}^{p})$.
Finally, recall the following map
\begin{align*}
   \xymatrix {\mathcal{X} \ar@{->}[r] & \mathbb{N}_{0}}; \qquad H\mapsto \mrdim_{\Bk} H\cap \Bk\mcG_{p}
\end{align*}
is upper semicontinuous \cite[Lemma 7.3]{Far1}. Thus  $\mathcal{Y}=\{H\in \mathcal{X} \mid \mrdim_{\Bk} H\cap \Bk\mcG_{p}\geq p^{\ell}-1 \}$ is closed.
\end{proof}

Suppose $\xymatrix{\iota_{\mcE}:\mcE \ar@{^{(}->}[r] & \mcG}$ is the canonical inclusion of an 
elementary abelian subgroup. Then there is a morphism 
$\xymatrix{\iota_{*,\mcE}: V_{r}(\mcE)\ar@{->}[r] & V_{r}(\mcG)}$ of support varieties. 
We set 
\begin{align*}
\mfC(\ell\uparrow,\mcG):=\bigcup_{r\geq \ell}\mfC(r,\mcG)
\end{align*}
and
\begin{align*}
V_{\mfC(\ell\uparrow,\mcG)}=\bigcup_{\mcE \in \mfC(\ell\uparrow,\mcG)} \iotaE .
\end{align*}

\begin{thm}{\label{P-V}}
Suppose $\mcG$ is an infinitesimal group scheme of height $\leq r$. Then
\begin{align*}
   \srk(\mcG)=\mmax\set{\ell \mid V_{r}(\mcG)=V_{\mfC(\ell\uparrow, \mcG)}}.
\end{align*}
\end{thm}
\begin{proof}
Let $\set{v_{0},\ldots,v_{p^{r}-1}}$ be the dual basis of the standard basis 
$\{T^{0},T^{1},\ldots,T^{p^{r}-1}\}$ of $\Bk[\mbG_{a(r)}]=\Bk[T]/(T^{p^{r}})$. Denote by $u_{i}=v_{p^{i}}$ for $0\leq i <r$,
this gives $\Bk\mbG_{a(r)}=\Bk[u_{0},\ldots,u_{r-1}]$. Now we turn to verify our statement. 
Proposition 3.8 of \cite{FP} gives a bijection 
\begin{align*}
  \xymatrix{\Theta_{\mathcal{G}}:\mathrm{Proj}(V_{r}(\mcG))\ar@{->}[r] & \mrP(\mcG)};\;\;[\alpha] 
  \mapsto 
                               [\alpha_{*}\circ \epsilon]
\end{align*}
where $\xymatrix{\epsilon: \Bk[u_{r-1}]\simeq \Ap \ar@{^{(}->}[r] & 
\Bk\mbG_{a(r)}}$, and $\xymatrix{\alpha_{*}: \Bk\mbG_{a(r)}\ar@{->}[r] & \Bk\mcG}$.
If $\alpha\in \mathrm{Pt}(\mcG)$, then it represents an equivalence class given by the map $\Theta_{\mcG}$, 
say $[\alpha]= \Theta_{\mcG}([\beta])$ for some $\beta\in V_{r}(\mcG)$. Further if we assume that 
$V_{r}(\mcG)=V_{\mfC(\ell\uparrow,\mcG)}$, then there exists a maximal elementary abelian subgroup
$\mcE\leq\mcG$ with $\cx_{\mcE}(\Bk)\geq \ell$ such that $\beta=\iota_{*,\mcE}(\gamma)$ for some 
$\gamma \in V_{r}(\mcE)$. Therefore, $[\alpha]=[\beta_{*}\circ \epsilon]=[\iota_{*,\Bk\mcE}\circ 
\gamma_{*}\circ\epsilon]= \iota_{*,\Bk\mcE}([\gamma_{*}\circ\epsilon])$ which lies in 
$\iota_{*,\Bk\mcE}(\mrP(\mcE))$, and this 
gives $\mrP(\mcG)=\mrP(\mcG)_{\ell}$. On the other hand, if $\alpha\in V_{r}(\mcG)\setminus\set{0}$ 
and suppose $\mrP(\mcG)=\mrP(\mcG)_{\ell}$, then $\Theta_{\mcG}([\alpha])$ is an equivalence class of 
$\mrP(\mcG)$. 
By our assumption there exists a maximal abelian unipotent subgroup $\mcU$ along with an 
elementary abelian
subgroup $\mcE_{\mcU}$ with $\cx_{\mcE_{\mcU}}(\Bk)\geq \ell$ and a $p$-point $\beta \in 
\mathrm{Pt}(\mcE_{\mcU})$ such that $\Theta_{\mcG}([\alpha])=\iota_{*,\Bk\mcE_{\mcU}}([\beta])$; see
\cite[Lemma 6.2.1]{Far4}. Again by the bijective map $\Theta_{\mcE_{\mcU}}$, we have 
$[\beta]=[\gamma_{*}\circ \epsilon]$ for some 
$\gamma\in V_{r}(\mcE_{\mcU})$, thus $\Theta_{\mcG}([\alpha])=\Theta_{\mcG}([\iota_{\mcE_{\mcU}}
\circ\gamma])$ and $\alpha=\iota_{\mcE_{\mcU}}\circ\gamma \in \iota_{*,\mcE_{\mcU}}(V_{r}
(\mcE_{\mcU}))$. Therefore, 
$V_{r}(\mcG)=V_{\mfC(\ell\uparrow,\mcG)}$, as desirable. 
\end{proof}

\begin{cor}
In Theorem \ref{P-V}, if $\mcG$ is of height $\leq 1$, then 
\begin{align*}
\srk(\mcG)=\mmax\set{\ell \mid V_{r}(\mcG)=V_{\mfC(\ell, \mcG)}}.
\end{align*}
\end{cor}
\begin{proof}
It suffices to show that any $\gamma\in V_{1}(\mcE)$ with $\cx_{\mcE}(\Bk) = s$ may factor through
an elementary subgroup $\mcE^{'}$ of $\mcE$ for $s^{'}\leq s$. Since $\mcG$ is of height $\leq 1$,
we have $\mcE\cong \mbG_{a(1)}^{\times s}$. The homomorphism $\xymatrix{\gamma: \mbG_{a(1)} 
\ar@{->}[r] & \mcE}$ is equivalent to $\xymatrix{d\gamma: \mfg_{a} \ar@{->}[r] & \mrlie(\mcE)}$
where $\mfg_{a}=\Bk \delta_{a}=\mrlie(\mbG_{a(1)})$. The image $d\gamma(\delta_{a})$ is contained in
an elementary subalgebra $\mfe_{s^{'}}$ of $\mrlie(\mcE)$ since $\mrlie(\mcE)$ is elementary abelian.
Thus, the image $\gamma(\mbG_{a(1)})$ is contained in an elementary group scheme $\underline{\mfe_{s^{'}}}$ of $\mcE$, as desirable.
\end{proof}

Suppose $\alpha\in V_{r}(\mcG)$. Consider the following set
\begin{align*}
\mfC(\ell\uparrow,\mcG)_{\alpha}:=\set{\mcE\in \mfC(\ell\uparrow,\mcG)\mid \alpha\in \iota_{*,\mcE}(V_{r}(\mcE))}
\end{align*}   
together with
\begin{align*}
r_{\alpha}^{\mcG}:=\max\set{\ell\mid \mfC(\ell\uparrow,\mcG)_{\alpha}\neq \emptyset}. 
\end{align*} 
Write $r_{min}^{\mcG}:=\min \set{r_{\alpha}^{\mcG}\mid \alpha\in V_{r}(\mcG)}$ and 
$\mcO_{rmin}^{\mcG}:=\{\alpha\in V_{r}(\mcG)\mid r_{\alpha}^{\mcG}=r_{min}^{\mcG}\}$. 

\begin{rem}{\label{auto}}
Let $\mcG=G_{r}$ be the Frobenius kernel of a smooth group scheme $G$. Then $G$ acts on 
$G_{r}$ via the adjoint representation. There results an action of $G$ on $V_{r}(G)$ and 
on $\mfC(\ell,G_{r})$.
Let $\alpha\in V_{r}(G)$ and $g\in G$. Based on these facts, we obtain
the relation $r_{\alpha}^{G_{r}}=r_{g.\alpha}^{G_{r}}$.
\end{rem}

\begin{lem}{\label{cloV}}
Suppose that $\mcG$ is an infinitesimal group schemes of height $\leq r$. Then
$V_{\mfC(\ell,\mcG)}$ is a closed subvariety of $V_{r}(\mcG)$. 
\end{lem}
\begin{proof}
We denote by $\mrpr_{1}$ the projection onto the first coordinate:
\begin{align*}
    \xymatrix{\mrpr_{1} : V_{r}(\mcG) \times \mfC(\ell,\mcG)\ar@{->}[r] & V_{r}(\mcG)}.
\end{align*}
Write $\Bk\mbG_{a(r)}=\Bk[u_{0},\ldots,u_{r-1}]$ as we did in Theorem \ref{P-V}.
Consider the set 
$\mcZ=\{(\alpha_{*},\Bk\mcE)\in Hom(\Bk\mbG_{a(r)}, \Bk\mcG)\times 
\mfC(p^{\ell}, \Bk\mcG) \mid \alpha_{*}(u_{i})\in \Bk\mcE, 0 \leq i < r\}$,
which is closed. Then by categorical equivalence the set
$\mcZ^{'}:=\{(\alpha,\mcE) \in V_{r}(\mcG)\times \mfC(\ell,\mcG) \mid \alpha \in 
\iota_{*,\mcG}(V_{r}(\mcE))\}$ is closed. 
Proposition \ref{proj} shows that $\mfC(\ell,\mcG)$ is complete.
Therefore, the image $V_{\mfC(\ell,\mcG)}=\mrpr_{1}(\mcZ^{'})$ of $\mcZ^{'}$
is closed in $V_{r}(\mcG)$ by general theory.
\end{proof}

\begin{thm}{\label{srk-inf}}
Suppose that $\mcG$ is an infinitesimal group scheme of height $\leq r$. Then $\srk(\mcG)=r_{min}^{\mcG}$ 
and $\mcO_{rmin}^{\mcG}$ is an open subset of $V_{r}(\mcG)$.
\end{thm}
\begin{proof}
Let $s=\srk(\mcG)$. Then $V_{r}(\mcG)=V_{\mfC(s\uparrow,\mcG)}$ and 
$\mfC(s\uparrow,\mcG)_{\alpha}\neq \emptyset$ for any $\alpha\in V_{r}(\mcG)$. 
Thus $r_{\alpha}^{\mcG}\geq s$, resulting $r_{min}^{\mcG}\geq s$. 
On the other hand, $r_{\alpha}^{\mcG}\geq r_{min}^{\mcG}$ 
gives $\mfC(r_{min}^{\mcG}\uparrow,\mcG)_{\alpha}\neq \emptyset$ for any $\alpha\in V_{r}(\mcG)$. 
Therefore, $V_{r}(\mcG)=V_{\mfC(r_{min}^{\mcG}\uparrow,\mcG)}$ and $s\geq r_{min}^{\mcG}$.

We now consider the function
\begin{align*}
\xymatrix{r^{\mcG}: V_{r}(\mcG) \ar@{->}[r] & \mathbb{N}};\; \alpha \mapsto \ r_{\alpha}^{\mcG}.
\end{align*}
Since $\mcB_{n}^{\mcG}:=\set{\alpha\in V_{r}(\mcG)\mid r_{\alpha}^{\mcG}\geq n}=V_{\mfC(n\uparrow,
\mcG)}=\bigcup_{s\geq n} V_{\mfC(s,\mcG)}$ is closed for every $n\in \mathbb{N}$; see Lemma \ref{cloV},
the function $r^{\mcG}$ is upper-semicontinuous(see \cite[Sect.1]{Far1}). 
As a result, $\mcO_{rmin}^{\mcG}=V_{r}(\mcG)\setminus \mcB_{\srk(\mcG)+1}^{\mcG}$ is open.
\end{proof}

\section{Infinitesimal group schemes: height $\leq 1$}
\subsection{}
{\label{general}}$\bf{General\;theory.}$
Let $\mcG$ be an infinitesimal group scheme of height $\leq 1$. Then $\mcG=\ug$ for some finite 
dimensional restricted Lie algebra $(\mfg,[p])$. Let $V(\mfg)$ be the fibre of zero of the map 
$\xymatrix{[p]:\mfg\ar@{->}[r] &\mfg}$, i.e.
\begin{align*}
   V(\mfg)=\set{x\in \mfg \; ;\; x^{[p]}=0}
\end{align*}
which is called the restricted nullcone of $\mfg$ and 
\begin{align*}
      V_{\mbE(r,\mfg)}:=\bigcup _{\mfe\in\mbE(r,\mfg)} \mfe \subseteq V(\mfg)
\end{align*}
where $\mbE(r,\mfg)$ is a projective variety defined in \cite{CFP}.
Accordingly, we define
\begin{align*}
 \mbE(r,\mfg)_{x}:=\set{\mfe \in \mbE(r,\mfg)  \; ;\; x \in \mfe},\;\;\;
 r_{x}^{\mfg}:=\max\set{r \; ;\;  \mbE(r,\mfg)_{x}\neq \emptyset}
\end{align*} 
and
\begin{align*}
   r_{min}^{\mfg}:= \min \set{ r_{x}^{\mfg} \; ;\; x\in V(\mfg)}.
\end{align*}

Being the special case of arbitrary infinitesimal group schemes, the passage from $\mrP(\ug)$ to 
$V(\mfg)$ is essentially a translation of the illustration as we did in section 2.3. 
We list the related results for $\mfg$:
\begin{itemize}
\item   $\srk(\ug)=\mmax \set{r \; ;\; V(\mfg)=V_{\mbE(r,\mfg)}}$.
\item   $\srk(\ug)=r_{min}^{\mfg}$. 
\item   $\mcO_{rmin}^{\mfg}:=\set{x\in V(\mfg) \; ; \; r_{x}^{\mfg}=r_{min}^{\mfg}}$ is open.
\end{itemize}

Without any real ambiguity, we denote the saturation rank  $\srk(\ug)$ of $\ug$ by $\srk(\mfg)$.

\begin{exmp}{\label{Heisenberg Lie algebra}}
We consider the $2n+1$-dimensional Heisenberg algebra $\mfh:=\bigoplus\limits_{i=1}^{n} \Bk x_{i}
\oplus\bigoplus\limits_{j=1}^{n}\Bk y_{j}\oplus \Bk z$, whose bracket and $p$-map are given by
\begin{align*}
    &x_{i}^{[p]}=0,\;y_{j}^{[p]}=0,\;z^{[p]}=0,\\
    &[x_{i},x_{j}]=0=[y_{i},y_{j}],\;[x_{i},y_{j}]=\delta_{ij}z,\;[z,\mathfrak{h}]=0
\end{align*}
respectively.
Suppose that $p \geq 3$, then Jacobson's formula implies that $\mfh$ is $[p]$-trivial, i.e. $V({\mfh})=
\mfh$. The proof of Proposition 2.2 in \cite{CFP} shows that the assignment 
$\phi: \mfe \mapsto \mfe/\Bk z$ is a bijection between the maximal abelian subalgebras of 
$\mfh$ and the maximal totally isotropic subspaces of the symplectic vector space $\mfh/\Bk z$.
Since every element of $\mfh/\Bk z$ is contained in a maximal totally isotropic subspace.
As a result, we have $\srk(\mfh)=n+1$.
\end{exmp}

\begin{rem}{\label{c-v}}
Let $\mfz_{\mfg}(x)=\set{y\in \mfg \; ;\; [y,x]=0}$ be the centralizer of $x$ in $\mfg$, and $\mfg_{x}^{\nil}:=V(\mfg)\cap \mfz_{\mfg}(x)$ be the intersection of $V(\mfg)$ and $\mfz_{\mfg}(x)$. For each positive integer $r$, we define
\begin{align*}
 \mathscr{C}_{r}(\mfg_{x}^{\nil}):=\set{(x_{1},x_{2},\ldots,x_{r})\in (\mfg_{x}^{\nil})^{\times r} \; ;\; [x_{i},x_{j}]=0 ,1\leq i,j\leq r} 
\end{align*}
to be the Zariski closed subvariety of $r$-tuples of pairwise commuting elements of $\mfg_{x}^{\nil}$, as well as 
$\mathscr{C}_{r}(\mfg_{x}^{\nil})^{\circ}$ the open subset of linear independent $r$-tuples of $\mathscr{C}_{r}(\mfg_{x}^{\nil})$. Denote by $r_{*}=\max\{r \; ;\; \mathscr{C}_{r}(\mfg_{x}^{\nil})^{\circ} \neq \emptyset \}$.
Then for every element $(x_{1},\ldots,x_{r_{*}})$ of $\mathscr{C}_{r_{*}}(\mfg_{x}^{\nil})$, 
$x\in \mrsp_{\Bk}\set{x_{1},\ldots,x_{r_{*}}}$. Otherwise, 
$(x,x_{1},\ldots,x_{r_{*}})\in \mathscr{C}_{r_{*}+1}(\mfg_{x}^{\nil})^{\circ}\neq \emptyset$, a contradiction.
Thus, each element in $\mathscr{C}_{r_{*}}(\mfg_{x}^{\nil})^{\circ}$ gives rise to an elementary subalgebra
which contains $x$. According to this, one can easily check that 
\begin{align*}
   r_{x}^{\mfg}=\max\set{r \; ;\; \mathscr{C}_{r}(\mfg_{x}^{\nil})^{\circ} \neq \emptyset },
\end{align*}
and there is a surjective map
\begin{align*}
   q:  \xymatrix{\mathscr{C}_{r_{x}^{\mfg}}(\mfg_{x}^{\nil})^{\circ} \ar@{->}[r] & \mbE(\mfg)_{x}}
\end{align*}
with the elements of $q^{-1}(\mfe_{x})$ for any $\mfe_{x}\in \mbE(\mfg)_{x}$ differ 
by the natural action of $\GL_{r_{x}}$.
\end{rem}

\subsection{}
{\label{reductive case}}
$\bf{Reductive\; Lie\; algebras.}$
In this section, we assume that $\mfg$ is the Lie algebra of a connected reductive algebraic 
$\Bk$-group $G$. Let $\msN(\mfg)$ be the nullcone of $\mfg$ consisting of all elements $x$ in 
$\mfg$ which are $[p]$-nilpotent in the sense that $x^{[p]^{r}}=0$ for some $r > 0$ depending 
on $x$. It is known that when $p$ is good for $G$, $\msN(\mfg)$ is finite union of $G$-orbits; 
see \cite[2.8.Theorem 1]{Jan2}.
The result also applies to $V(\mfg)$, by the fact that $V(\mfg)$ is a $G$-stable subvariety of 
$\msN(\mfg)$. We list the good primes
and their Coxeter number for $G$ being simple algebraic groups in the following Table 
\ref{Good-Coxe}. With the finiteness property behind us, the existence of three canonical nilpotent 
orbits should be known: the regular (or principal) orbit $\mcO_{reg}$ (or $\mcO_{prin}$), the 
subregular orbit $\mcO_{subreg}$ and the minimal orbit $\mcO_{min}$; see \cite[4.1-4.3]{CM}
\footnote{We often refer to \cite{CM} in this Section. Although the results there are obtained
over complex numbers, they are still valid in positive characteristic $p$ as long as $p$ is good. 
See \cite[Sect 3.9]{CLNP}.}.
A connected reductive group $G$ is said to be \emph{standard} if it 
satisfies the follwoing hypotheses(see \cite[2.9]{Jan2}):
\begin{itemize}
\item  The derived subgroup $G^{(1)}$ of $G$ is simply-connected
\item  $p$ is a good prime for G.
\item  The Lie algebra $\mfg$ of $G$ has a non-degenerate symmetric bilinear $G$-invariant form 
           $\xymatrix{B:\mfg\times\mfg\ar@{->}[r] &\Bk}$.
\end{itemize}

\begin{table}[ht]
\centering
\begin{tabular}[b]{|r||c|c|c|c|c|c|c|}
\hline
Type & $A_{n}$ & $B_{n}/C_{n}$
& $D_{n}$  & $E_{6}/F_{4}$ & $E_{7}$ & $E_{8}$ & $G_{2}$ \\
\hline
Good primes $p$ & $\geq 2$ & $\geq 3$ & $\geq 3$ & $\geq 5$ & $\geq 5$ & $\geq 7$ &
                $\geq 5$ \\
\hline
Coxeter number $h(G)$ &$n+1$ & $2n$ & $2n-2$ & $12$ & $18$ & $30$ & $6$ \\
\hline
\end{tabular}
\caption{Good primes and Coxeter  number}
\label{Good-Coxe}
\end{table}
Recall that the semisimple rank $\mrrk(G)$ of $G$ is defined to be the rank of its derived subgroup $G^{(1)}$.
The main theorem of this section is the following:
\begin{thm}{\label{ss-rank}}
Let $G$ be a standard connected reductive algebraic group with $\mfg:=\mrlie(G)$. Assume that $p\geq h$, where $h$ is the Coxeter number of $G$. Then 
\begin{align*}
 \srk(\mfg)=\mrrk(G).
\end{align*}
\end{thm}
\begin{proof}
Let $G^{(1)}=[G,G]$ be the derived subgroup of $G$ and $\mfg^{(1)}$ be the corresponding Lie algebra. By our assumption on $G$, $G^{(1)}$ is simply-connected and semisimple. Let $G_{1}, G_{2},\ldots,G_{m}$ be the simple simply-connected normal subgroups of $G^{(1)}$ with $\mfg_{i}=\mrlie(G_{i})$. Then $G^{(1)}=G_{1}\times G_{2}\times\cdots \times G_{m}$ as well as $\mfg^{(1)}=\mfg_{1}\oplus\mfg_{2}\oplus\cdots\oplus \mfg_{m}$ \cite[Sect.3.1]{Pre}. 
Moreover, $G_{i}$ is one of the cases: 
\begin{enumerate}
 \item $G_{i}$ is simple, simply-connected and not of type $A_{\Bk p-1}$ ;
 \item $G_{i}=\SL(V_{i})$ and $p \mid \mrdim V_{i}$.
\end{enumerate}
By case (ii), in accordance with the "standard" criterion, we define groups $G_{i}^{'}$ by setting
\begin{align*}
   G_{i}^{'}=
      \begin{cases}
         \GL(V_{i}),  &\;\text{if}\;\;G_{i}=\SL(V_{i}) \;\mbox{and}\; p \mid \mrdim V_{i}, \\
         G_{i},  &\;\mbox{otherwise}. 
      \end{cases}
\end{align*}
and further denote by $G^{'}=G_{1}^{'}\times G_{2}^{'}\times \cdots \times G_{m}^{'}$. 
In this fashion, one can check that $G^{'}$ satisfies the three hypotheses, that is, it is standard. 
Let $\mfg_{i}^{'}=\mrlie(G_{i}^{'})$, and $\mfg^{'}=\mrlie(G^{'})$, then we have
$\mfg^{'}=\mfg_{1}^{'}\oplus\mfg_{2}^{'}\oplus\cdots\oplus \mfg_{m}^{'}$. According to 
\cite[6.2]{GP}, there are tori $T_{0}$ and $T_{1}$ with their corresponding Lie algebras $\mft_{0}$ 
and $\mft_{1}$ such that $\mfg^{'}\oplus \mft_{0}=\mfg\oplus \mft_{1}$. 
By exploiting this result, we have  $V(\mfg)=V(\mfg^{'})$ as well as 
$\mbE(r,\mfg)=\mbE(r,\mfg^{'})$ for any positive integer $r$, which 
implies $\srk(\mfg)=\srk(\mfg^{'})$.

We are in a position to compute $\srk(\mfg^{'})$. Notice that $G$ and $G^{'}$ have the same Coxeter 
number. If $p\geq h$, then $V(\mfg^{'})=\msN(\mfg^{'})=\overline{\mcO_{reg}}$ 
is irreducible ensured by \cite[Lemma 6.2]{Jan2}, where 
$\mcO_{reg}=G^{'}.e_{reg}$ and $e_{reg}$ is a regular nilpotent element of $\mfg^{'}$. 
By Section \ref{general}, we have
$\mcO_{reg} \cap \mcO_{rmin}^{\mfg^{'}} \neq \emptyset$. 
By Remark \ref{auto}, there is $r_{x}^{\mfg^{'}}=r_{h.x}^{\mfg^{'}}$ for any 
$h\in G^{'}$ and $x\in \mfg^{'}$. 
As a result, $\mcO_{reg}\subseteq \mcO_{rmin}^{\mfg^{'}}$ and
$\srk(\mfg^{'})=r_{min}^{\mfg^{'}}=r_{e_{reg}}^{\mfg^{'}}$.

Write $e_{reg}=e_{1}+ e_{2}+\cdots + e_{m}, e_{i}\in \mfg_{i}^{'}$. Then
\begin{align*}
C_{G^{'}}(e_{reg})=C_{G_{1}^{'}}(e_{1})\times C_{G_{2}^{'}}(e_{2})\times \cdots 
 \times C_{G_{m}^{'}}(e_{m}) 
\end{align*}
and 
\begin{align*}
 \mfz_{\mfg^{'}}(e_{reg})=\mfz_{\mfg_{1}^{'}}(e_{1})\oplus \mfz_{\mfg_{2}^{'}}(e_{2})
 \oplus \cdots \oplus \mfz_{\mfg_{m}^{'}}(e_{m})
\end{align*}
One can check that $e_{reg}$ is regular nilpotent in $\mfg^{'}$ if and only if
each $e_{i}$ is regular nilpotent in $\mfg_{i}^{'}$.
Since each $G_{i}^{'}$ is standard, there is 
$\mfz_{\mfg_{i}^{'}}(e_{i})=\mrlie\;C_{G_{i}^{'}}(e_{i})$. 
So if each $C_{G_{i}^{'}}(e_{i})$ is abelian, then $\mfz_{\mfg^{'}}(e_{reg})$ is abelian.
Check that each $G_{i}^{'}$ is $D$-standard reductive group 
(See \cite[Definition 3.2]{McT}, \cite[Remark 3]{McN} and \cite[Remark 2.5.6(a)]{Let}), the results of
\cite[(5.2.4)]{McT} implies the abelian property of $\mfz_{\mfg^{'}}(e_{reg})$.
Due to the Jordan-Chevalley decomposition, there is a direct sum of Lie algebra
$\mfz_{\mfg^{'}}(e_{reg})=(\mfz_{\mfg^{'}}(e_{reg}))_{s}\oplus(\mfz_{\mfg^{'}}(e_{reg}))_{n}$. 
Since regular nilpotent elements are distinguished(see \cite[Sect. 7.13]{Hum}), this gives rise to
$(\mfz_{\mfg^{'}}(e_{reg}))_{s}=\mfz(\mfg^{'})$ by \cite[Theorem 3]{Lev}. At this point, the abelian
property of $(\mfz_{\mfg^{'}}(e_{reg}))_{n}$ ensures that, it is the maximal elementary subalgebra of 
$\mfg^{'}$ containing $e_{reg}$ of dimension $(\dim \mfz_{\mfg^{'}}(e_{reg})-\dim \mfz(\mfg^{'}))$.
For $G^{'}$ and $\mfg^{'}$, we have $\dim \mfz_{\mfg^{'}}(e_{reg})=\dim C_{G^{'}}(e_{reg})=\rk(G^{'})$
since $G^{'}$ is standard.
That gives, $\srk(\mfg^{'})=\dim (\mfz_{\mfg^{'}}(e_{reg}))_{n}=\rk(G^{'})-\dim \mfz(\mfg^{'})$.
  
We compute $\rk(G^{'})$ and $\dim \mfz(\mfg^{'})$.
If $G_{i}^{'}=G_{i}$, then $\mfg_{i}^{'}$ is simple and $\mfz(\mfg_{i}^{'})=0$. 
If $G_{i}^{'}=\GL(V_{i})$ with $p\mid\mrdim V_{i}$, then $\mfg_{i}^{'}$ has a one-dimensional center.
By denoting
\begin{align*}
rk_{p}(G):=\ord\set{G_{i} \; ;\; G_{i}=\SL(V_{i})\;\mbox{and}\; p\mid \mrdim V_{i}}
\end{align*}
we have $\dim \mfz(\mfg^{'})=rk_{p}(G)$. 
Since $\rk(\GL(V_{i}))=1+\rk(\SL(V_{i}))$, this implies
$\rk(G^{'})=\rk(G^{(1)})+rk_{p}(G)$. 
As a result,
\begin{align*}
\srk(\mfg^{'})=\rk(G^{'})- rk_{p}(G)
                = (\rk(G^{(1)})+ rk_{p}(G))- rk_{p}(G) = \mrrk(G),                           
\end{align*}
that is our $\srk(\mfg)$ and we finish the proof.
\end{proof}

\subsection{}
$\bf{Open\;set\;\mcO_{rmin}}$.
We specialise our $G$ in Section \ref{reductive case} into $\SL_{n}(\Bk)$ and $n\geq 3$. 
Although $\SL_{n}(\Bk)$ is not standard when $p \mid n$, 
we find from the proof of Theorem \ref{ss-rank} that $\srk(\mfsl_{n}(\Bk))=n-1$ whenever $p\geq n$. 
We now recollect some material on nilpotent orbits from the partition point of view 
for $G$. Let $\mcP(n)$ be a set of partitions of $n$, with
which we can endow a dominance partial order $\trianglelefteq$.
Given two partitions $\lambda$ and $\mu$, we say $\lambda$ dominates $\mu$, 
provided $\mu \trianglelefteq \lambda$.
If $\lambda=(\lambda_{1},\ldots,\lambda_{t})$ is such a partition, then we assign to it a nilpotent 
matrix $x_{\lambda}=\mathrm{diag}(N_{1},\ldots,N_{t})$ with upper triangular Jordan blocks 
$N_{1},\ldots,N_{t}$ with sizes $\lambda_{1}\times\lambda_{1},\ldots,\lambda_{t}\times \lambda_{t}$, 
and put $\mcO_{\lambda}=G.x_{\lambda}$ as a nilpotent $G$-orbit.
It is known that the nilpotent orbits of $G$ in $\mfg$ can be
described in terms of partitions, and if $\mcO_{\lambda}$ and $\mcO_{\mu}$ are two nilpotent 
orbits in $\mfg$, then $\mcO_{\mu}\subseteq \overline{\mcO_{\lambda}}$ 
if and only if $\mu \trianglelefteq \lambda$; see \cite[Theorem 5.1.1/6.2.5]{CM}.
Write $n=qp+r$ with $0\leq r <p$. Denote by
$\lambda \vdash n$ the partition with $q$ parts of size $p$ and 1 part of size $r$. Then 
$\lambda$ is maximal with respect to $\trianglelefteq$ among any other partitions of $n$.
As a result, $V(\mfg)=\overline{\mcO_{\lambda}}$.

Inspired by Remark \ref{c-v} and the proof of Theorem \ref{ss-rank}, the determination 
of the centralizer of a nilpotent element plays a central role in our computation.
By the natural embedding of $\mfsl_{n}(\Bk)$ in $\gl(\mbV)$ with $\mrdim \mbV =n$, 
we begin with an observation from $\gl(\mbV)$. Let $e\in \gl(\mbV)$ be nilpotent with the corresponding partition
$(\lambda_{1},\ldots,\lambda_{t})$, $\mfz(e)$ be the centralizer of $e$ in $\gl(\mbV)$. 
It is assumed that $\lambda_{1}\geq \lambda_{2}\geq \cdots
\geq \lambda_{t}>0$. Then there exist elements $v_{1},\ldots,v_{t}\in \mbV$ such that all 
$e^{j}.v_{i}$ with $1\leq i \leq t$ and $0\leq j <\lambda_{i}$ form a basis of $\mbV$ together with
$e^{\lambda_{i}}.v_{i}=0$. Let $\xi\in\mfz(e)$, then $\xi$ is completely determined by $\xi(v_{i})$ with
$1\leq i\leq t$ because $\xi(e^{j}.v_{i})=e^{j}.\xi(v_{i})$, but $e^{\lambda_{i}}.\xi(v_{i})=0$.
One can easily check that 
\begin{align*}
 \xi(v_{i})=\sum_{j=1}^{t}\sum_{s=\max\set{\lambda_{j}-\lambda_{i},0}}^{\lambda_{j}-1}a_{ijs}e^{s}.v_{j}
\end{align*}
for some $a_{ijs}$, which gives the basis $\set{\xi_{i}^{j,s}}$ of $\mfz(e)$ defined by
$$
\left\{
\begin{array}{l}
\xi_i^{j,s}(v_i)=e^s{\cdot}v_j, \\
\xi_i^{j,s}(v_r)=0 \enskip \mbox{for } r\neq i, \\
\end{array}\right.
\quad 1\leq i,j\le k, \ \mbox{ and }\ \max\{\lambda_j-\lambda_i, 0\} \leq s< \lambda_j \ .
$$
It is convenient to assume that $\xi_{i}^{j,s}=0$ whenever $s$ is not within the appropriate bound.

Given two basis elements $\xi_{i}^{j,s}$ and $\xi_{p}^{q,r}$ of $\mfz(e)$, we have the composition rule:
$\xi_{i}^{j,s}\cdot\xi_{p}^{q,r}=\delta_{q,i}\xi_{p}^{j,s+r}$, and further their bracket:
\begin{align*}
   [\xi_{i}^{j,s},\xi_{p}^{q,r}]=\delta_{q,i}\xi_{p}^{j,s+r}-\delta_{j,p}\xi_{i}^{q,s+r}.
\end{align*}

\begin{lem}{\label{subreg}}
Suppose $G=\SL_{n}(\Bk)(n\ge 3)$. Let $\tau=(n-1,1)\vdash n$ be the partition corresponding to the 
subregular nilpotent orbit of $\mfg$, $x_{\tau}$ be the nilpotent matrix given by the partition $\tau$. 
If $p\geq n-1$, then we have $r_{x_{\tau}}=n-1$.
\end{lem}
\begin{proof} 
Let $e:=x_{\tau}$ with $\tau=(n-1,1)$, a nilpotent element in $\gl(\mbV)$ in a natural way. Keep the 
notation for $\mfz(e)$ as above, the basis of $\mfz(e)$ as described is : 
\begin{align*}
   \xi_{1}^{1,s}, 0\leq s\leq n-2 ;  \;\;\; \xi_{1}^{2,0};\;\;\;\xi_{2}^{1,n-2};\;\;\;\xi_{2}^{2,0}.
\end{align*}
Let $\mfz(e)^{'}=\mfz(e)\cap \mfg$, the centralizer of the nilpotent element $e$ in $\mfg$, and
$\xi\in\mfz(e)^{'}$. Write 
\begin{align*}
\xi=\sum_{s=0}^{n-2}a_{11s}\xi_{1}^{1,s} + a_{120}\xi_{1}^{2,0}+
a_{21(n-2)}\xi_{2}^{1,n-2}+a_{220}\xi_{2}^{2,0}.
\end{align*}
Since it is an element of $\mfg$, this implies $a_{220}=-(n-1)a_{110}$. Thus the basis of 
$\mfz(e)^{'}$ is 
\begin{align*}
   \xi_{1}^{1,s}, 1\leq s\leq n-2;\;\;\; \xi_{1}^{2,0};\;\;\;\xi_{2}^{1,n-2};\;\;\;\xi_{1}^{1,0}-(n-1)\xi_{2}^{2,0}.
\end{align*}
Keep the form of $\xi$ with the relation of $a_{110}$ and $a_{220}$. We want to determine the nilpotent part
$\mfz(e)^{'}_{n}$ of $\mfz(e)^{'}$. Given by the composition rule, there are the relations:
(1)\;$\xi_{1}^{1,0}\cdot\xi_{1}^{1,0}=\xi_{1}^{1,0}$, $\xi_{2}^{2,0}\cdot\xi_{2}^{2,0}=\xi_{2}^{2,0}$;
(2)\;$\xi_{2}^{1,n-2}\cdot\xi_{1}^{2,0}=\xi_{1}^{1,n-2}$;
(3)\;$(\xi_{1}^{1,s})^{n-1}=0$ for $1\leq s\leq n-2$, and $(\xi_{1}^{2,0})^{2}=(\xi_{2}^{1,n-2})^{2}=0$.
Depending on these rules, we conclude that if
$n>3$ and $p\geq n-1$, or $n=3$ and $p> n-1$, then $\mfz(e)^{'}_{n}$ is a vector space having the following basis
\begin{align*}
  \xi_{1}^{1,s}, 1\leq s\leq n-2 ;\;\;\; \xi_{1}^{2,0};\;\;\;\xi_{2}^{1,n-2}.
\end{align*}
Alternatively, in the case $n=3$ and $p=2$, $\xi\in\mfz(e)^{'}_{n}$ with $\xi^{p^{2}}=0$ if and only if $a_{110}=a_{220}=0$, with $\xi^{p}=0$ if and only if $a_{110}=a_{220}=a_{120}\cdot a_{211}=0$.
Now we are able to determine the maximal elementary subalgebras $\mfe_{\tau}$ that will be assigned to $e$.
Observe that for $1\leq s \leq n-2$, 
$[\xi_{1}^{1,s},\xi_{1}^{1,s^{'}}]=\xi_{1}^{1,s+s^{'}}-\xi_{1}^{1,s+s^{'}}=0$ for $1\leq s^{'}\leq n-2$, $[\xi_{1}^{1,s},\xi_{1}^{2,0}]=-\xi_{1}^{2,s}=0$ and 
$[\xi_{1}^{1,s},\xi_{2}^{1,n-2}]=\xi_{2}^{1,s+n-2}=0$. 
We then have $[\xi_{1}^{1,s},\xi]=0$ for any $\xi\in\mfz(e)^{'}_{n}$ and $1\leq s\leq n-2$.
Still, note that $[\xi_{1}^{2,0},\xi_{2}^{1,n-2}]\neq 0$. 
Thus, if $n>3$ and $p\geq n-1$, or $n=3$ and $p>n-1$, then 
\begin{align*}
\mfe_{\tau}=\Bk \xi_{1}^{1,1}\oplus\cdots\oplus\Bk\xi_{1}^{1,n-2}\oplus \Bk(a\xi_{1}^{2,0}+b\xi_{2}^{1,n-2});
\end{align*}
parameterized by $(a:b)\in \mbP^{1}$ and if $n=3$ and $p=2$, then
\begin{align*}
  \mfe_{\tau}=\Bk \xi_{1}^{1,1}\oplus\cdots\oplus\Bk\xi_{1}^{1,n-2}\oplus \Bk\xi^{*}, \;\mbox{where}\; 
  \xi^{*}\in\set{\xi_{1}^{2,0},\xi_{2}^{1,n-2}}
\end{align*}
Notice that all of them show $r_{e}=r_{x_{\tau}}=n-1$, and this completes our proof.
\end{proof}

\begin{rem}
We remark here that if $p=n-1$, then $V(\mfg)=\overline{\mcO_{\tau}}$ for $\tau=(n-1,1)$. According to 
the procedure processed in Theorem \ref{ss-rank} and the result of Lemma \ref{subreg} we have
$\srk(\mfg)=r_{x_{\tau}}=n-1$.
This shows that the equality $\srk(\mfg)=\mrrk(G)$ still holds for smaller $p$, like $p=n-1$.
\end{rem}

\begin{lem}{\label{fuzzy}}
Suppose $G=\SL_{n}(\Bk)$ with $n\geq 3$. Let $\mcO_{\lambda}\subseteq \msN(\mfg)\setminus 
(\mcO_{reg}\sqcup\mcO_{subreg})$ be a nilpotent orbit, and $x_{\lambda}\in \mcO_{\lambda}$ 
be the corresponding nilpotent element given by the partition $\lambda$. 
If $p\geq \mmax\set{2, n-2}$, then we have $r_{x_{\lambda}}\geq n$.
Specifically, when $\lambda=(2,1^{n-2})$ or $\lambda=(1^{n})$ we have 
$r_{x_{\lambda}}=\lfloor \frac{n^{2}}{4}\rfloor$.
\end{lem}
\begin{proof}
According to the dominance order $\unlhd$, we know that $\mcO_{\lambda}\subseteq \msN(\mfg)\setminus 
(\mcO_{reg}\sqcup\mcO_{subreg})$ if and only if $\lambda\unlhd (n-2,2)$. 
Write $\lambda=(\lambda_{1},\ldots,\lambda_{t})$ with $\lambda_{1}\geq \lambda_{2}\geq 
\cdots\geq\lambda_{t}>0$.  We are going to read off $x_{\lambda}$ from the basis $\set{\xi_{i}^{j,s}}$ 
of $\gl(\mbV)$. Let $s=\max\set{i \; ;\; \lambda_{i}\geq 2}$. If $1\leq i\leq s$, the action of
$\xi_{i}^{i,1}$ on $\mbV$ gives an unique non-zero Jordan block of size $\lambda_{i}\times\lambda_{i}$.
Then we can immediately know
\begin{align*}
  x_{\lambda}=\xi_{1}^{1,1}+\xi_{2}^{2,1}+\cdots+\xi_{s}^{s,1}.
\end{align*}
In what follows, we consider three cases to estimate the local saturation rank of $x_{\lambda}$.
In case of all $\lambda_{i}\geq 2$, we find there is an elementary subalgebra
\begin{align*}
\bigoplus_{i=1}^{t}\bigoplus_{r=1}^{\lambda_{i}-1}\Bk \xi_{i}^{i,r}\oplus\bigoplus_{i=1}^{t-1}
   \Bk\xi_{i}^{i+1,\lambda_{i+1}-1}\oplus \Bk\xi_{t}^{1,\lambda_{1}-1}
\end{align*}
that contains $x_{\lambda}$ with dimension $\sum_{i=1}^{t}\lambda_{i}=n$.
Alternatively, there is at least one block with size $1\times 1$. If there is only one, then
$s=t-1$ and $\lambda=(\lambda_{1},\ldots,\lambda_{s},1)$. Since $\lambda\unlhd (n-2,2)$, then $s\geq 2$ and
we find an elementary subalgebra
\begin{align*}
\bigoplus_{i=1}^{s}\bigoplus_{r=1}^{\lambda_{i}-1}\Bk \xi_{i}^{i,r}\oplus\bigoplus_{i=1}^{s-1}
\Bk\xi_{i}^{i+1,\lambda_{i+1}-1}\oplus\Bk\xi_{t}^{s,\lambda_{s}-1}\oplus\Bk\xi_{t}^{1,\lambda_{1}-1}
\end{align*}
which contains $x_{\lambda}$. 
Otherwise, $\lambda$ has the form $(\lambda_{1},\ldots,\lambda_{s},1^{t-s})$
for $t-s\geq 2$. Let $\xi=\sum_{i=1}^{t-s-1}\xi_{s+i+1}^{s+i,0}$.
Then we assign to $x_{\lambda}$ an elementary subalgebra
\begin{align*}
\bigoplus_{i=1}^{s}\bigoplus_{r=1}^{\lambda_{i}-1}\Bk \xi_{i}^{i,r}\oplus\bigoplus_{i=1}^{t-s-1}\Bk \xi^{i}
\oplus\bigoplus_{i=1}^{s}\Bk\xi_{i}^{i+1,\lambda_{i+1}-1}
\oplus \Bk\xi_{t}^{1,\lambda_{1}-1}.
\end{align*}
where $\xi^{i}$ is the $i$-th power of $\xi$.
One can compute that $r_{x_{\lambda}}\geq n$ for both cases.

We now consider $\lambda=(1^{n})$ or $\lambda=(2,1^{n-2})$.
When $\lambda=(2,1^{n-2})$, $\mcO_{\lambda}$ is the $G$-orbit of the highest root $E_{1n}$. 
When $n=2m$(resp. $n=2m+1$), Theorem 2.7(resp. Theorem 2.8) in \cite{CFP} shows that the 
elementary subalgebra which contains $E_{1n}$ has maximal dimension $m^{2}$(resp. $m(m+1)$). 
Therefore, we have $r_{x_{\lambda}}=\lfloor \frac{n^{2}}{4}\rfloor$.
\end{proof}

\begin{rem}
Suppose that $n\geq 4$. If $p=n-2$, then $V(\mfg)=\overline{\mathcal{O}_{\lambda}}$ with 
$\lambda=(n-2,2)$.
By the same token, we have $\srk(\mfg)=r_{x_{\lambda}}\geq n$ according to Lemma \ref{fuzzy}. 
This shows that for smaller $p<n$, like $p=n-2$, we have $\srk(\mfg)> \mrrk(G)$.
\end{rem}

\begin{thm}
Suppose $G=\SL_{n}(\Bk)$ with $n\geq 3$. If $p\geq n$, then
\begin{align*}
 \mcO_{rmin}=\mcO_{reg}\bigsqcup\mcO_{subreg}.
\end{align*}
\end{thm}
\begin{proof}
Theorem \ref{ss-rank} tells us $\mcO_{rmin}=\set{x\in V(\mfg)\; ;\; r_{x}=n-1}$. Consecutive applications
of Lemma \ref{subreg} and Lemma \ref{fuzzy} yield $\mcO_{rmin}=\mcO_{reg}\bigsqcup\mcO_{subreg}$.
\end{proof}


\section{Infinitesimal group schemes: $\SL_{n(2)}$}
\subsection{}
$\bf{Inequality}.$ 
In Section 3 we described the behavior of restricted Lie algebras, 
emphasising on reductive Lie algebras. So in this Section
we will first show the relation between the infinitesimal group schemes $\mcG$
of height $\leq r$ and their first Frobenius kernel $\mcG_{(1)}$ in their respective
saturation ranks.

\begin{lem} Let $\mcE$ be an elementary abelian group scheme. Then there exist 
$r_1,\ldots,r_n$ and $r \in \mathbb{N}$ such that
\[ \mcE \cong \prod_{i=1}^n \mbG_{a(r_i)}\!\times\!E_r.\]
\end{lem}

\begin{proof} Let $\mathfrak{D}(\mcE)$ be the Cartier dual of the abelian group scheme $\mcE$, so that 
$\Bk[\mathfrak{D}(\mcE)] \cong \Bk\mcE$. 
It follows that $\mathfrak{D}(\mcE)$ is an infinitesimal group of height $1$, implying
that its Frobenius morphism $\mathfrak{F}_{\mathfrak{D}(\mcE)}$ is zero. 
Letting $\mathfrak{V}_\mcE$ be the Verschiebung of $\mcE$ (See\ \cite[(IV,\S3,${\rm n}^{\rm o}4$)]{DG}), we obtain, observing \cite[(IV,\S3,4.9)]{DG},
\[ \mathfrak{D}(\mathfrak{V}_\mcE) = \mathfrak{F}_{\mathfrak{D}(\mcE)} =0.\]
Consequently, $\mathfrak{V}_\mcE=0$ and \cite[(IV,\S3,6.11)]{DG} yields the asserted isomorphism. \end{proof}

\begin{lem}{\label{inequa}}
Let $\mcG$ be an infinitesimal group scheme of height $\leq r$. Then 
$\srk(\mcG)\leq r \cdot \srk(\mcG_{(1)})$.
\end{lem}
\begin{proof}
Let $\mcE\subseteq \mcG$ be an elementary abelian subgroup. Then $ht(\mcE)\leq r$, so that
$\mcE=\prod_{i=1}^{r} \ell_{i} \mbG_{a(i)}$. Thus, 
$\cx_{\mcE}(\Bk)=\sum_{i=1}^{r}\ell_{i}\cdot i \leq r\cdot(\sum_{i=1}^{r}\ell_{i})$.
Moreover, $\mcE_{(1)}=\prod_{i=1}^{r}\ell_{i}\mbG_{a(1)}$, so that $\cx_{\mcE_{(1)}}
(\Bk)=\sum_{i=1}^{r}\ell_{i}$. This gives $\cx_{\mcE}(\Bk)\leq r\cdot \cx_{\mcE_{(1)}}(\Bk)$.

Assume that $\srk(\mcG)> r\cdot \srk(\mcG_{(1)})$. Then 
\begin{align*}
  V_{r}(\mcG)= \bigcup_{\cx_{\mcE}(\Bk)> r\cdot \srk(\mcG_{(1)})} V_{r}(\mcE) \;\;\mbox{and} \;\; 
  V_{1}(\mcG_{(1)})=\bigcup_{\cx_{\mcE}(\Bk)> r\cdot \srk(\mcG_{(1)})} V_{1}(\mcE_{(1)})
\end{align*}
This gives $\srk(\mcG_{1})< \cx_{\mcE_{(1)}}(\Bk)$, a contradiction.
\end{proof}

\subsection{}
$\bf{Equality}.$
We now suppose that $\mcG=\SL_{n(2)} \equiv \Ker\{F^{2}:\xymatrix{\SL_{n}\ar@{->}[r] &\SL_{n}\}}$,
where the geometric Frobenius $F:\SL_{n}(R)\rightarrow\SL_{n}(R)$ is 
defined by raising  each matrix entry to the $p^{th}$ power. 
Let $\msN=\msN(\mfsl_{n}(\Bk))$ be the nilpotent variety of $\mfsl_{n}(\Bk)$
and we consider the nilpotent commuting variety(see \cite{Pre}). 
\begin{align*}
\cnil(\mfsl_{n}(\Bk)):=\set{(x,y)\in \msN\times \msN \mid [x,y]=0}.
\end{align*}

The variety of 1-parameter subgroup $V_{2}(\SL_{n(2)})$ of $\SL_{n(2)}$ has been 
detected in terms of 2-tuples of pairwise commuting $p$-nilpotent matrices; see 
\cite[Proposition 1.2/Lemma 1.8]{SFB1}. 
More specifically, if $p\geq n$ then the element of $V_{2}(\SL_{n(2)})$ has uniquely
the form 
\begin{align*}
\xymatrix{\exp_{\underline{\alpha}}: \mbG_{a(2)}\ar@{->}[r] &\SL_{n(2)}} 
\end{align*}
where $\alpha=(\alpha_{0},\alpha_{1})\in \cnil(\mfsl_{n}(\Bk))$ by sending for 
any $\Bk$-algebra $A$ and $s\in A$
to $\exp_{\underline{\alpha}}(s)=\exp(s\alpha_{0})\cdot \exp(s^{p}\alpha_{1})$. 
Here for any $p$-nilpotent matrix $x\in \mfsl_{n}(A)$, we set
\begin{align*}
  \exp(x)=1+x+\frac{x^{2}}{2}+\cdots + \frac{x^{p-1}}{(p-1)!}\in \SL_{n}(A).
\end{align*}

Let $e\in \msN$ be a nilpotent element in $\mfsl_{n}(\Bk)$. 
As in section 3.2, we denote by $\mfz(e):=\mfz_{\mfsl_{n}(\Bk)}(e)$ the centralizer 
of $e$ in $\mfsl_{n}(\Bk)$, by $\mfC(e)$ the Zariski closure of 
$\SL_{n}(\Bk).(e,\msN\cap \mfz(e))$ and by $\mcO_{e}$ the 
$\SL_{n}(\Bk)$-orbit of $e$. 
Consider the morphism 
\begin{align*}
  \xymatrix{ \xi: \SL_{n}(\Bk) \times (\msN\cap \mfz(e)) \ar@{->}[r] & 
                              \mfC(e); \;\; \xi(g,x)=(\Ad(g).e,\Ad(g).x)}
\end{align*}
which is dominant, the canonical embedding 
\begin{align*}
\xymatrix{\iota:\msN\cap \mfz(e) \ar@{^{(}->}[r] &
\SL_{n}(\Bk)\times (\msN\cap \mfz(e))}
\end{align*}
which maps $x$ to $\iota(x)=(1,x)$ and their composition
$\xymatrix{\xi\circ\iota: \msN\cap \mfz(e) \ar@{->}[r] & \mfC(e)}$. 
Now we assume that $p\geq n$ and $e$ is a regular nilpotent element of $\mfsl_{n}(\Bk)$. 
Then $\mfz(e)=\mrsp_{\Bk}\set{e,e^{2},\ldots,e^{n-1}}$, implying
$\msN \cap \mfz(e)=\mbA^{n-1}$ is irreducible.

\begin{lem}{\label{non-empty}}
$(\xi\circ \iota)^{-1}(\mcO_{rmin}^{\SL_{n(2)}})$ is an non-empty open subset of $\msN \cap \mfz(e)$.
\end{lem}
\begin{proof}
We first prove the emptyness.
Since $\xi$ is a dominant morphism of irreducible affine varieties, the image
$\SL_{n(2)}(\Bk).(e,\msN\cap \mfz(e))$ contains a non-empty open subset $U$ of $\mfC(e)$.
Then the non-trivial intersection $\mcO_{rmin}^{\SL_{n(2)}}\cap U$ implies that 
$\mcO_{rmin}^{\SL_{n(2)}}$ contains an element $g^{'}.(e,x^{'})$ for some 
$g^{'}\in \SL_{n(2)}(\Bk)$ and $x^{'}\in \msN\cap \mfz(e)$. 
Notice that $\exp_{\underline{g^{'}.(e,x^{'})}}=g^{'}.\exp_{\underline{(e,x^{'})}}$, 
Remark \ref{auto} ensures that $(e,x^{'})\in \mcO_{rmin}^{\SL_{n(2)}}$. 
As a result, $x^{'}\in (\xi\circ \iota)^{-1}(\mcO_{rmin}^{\SL_{n(2)}})$ and 
$(\xi\circ \iota)^{-1}(\mcO_{rmin}^{\SL_{n(2)}})$ in not empty, as desired.

Then we show it is open. Note that $\mcO_{rmin}^{\SL_{n(2)}}$ is an open subset of $\mfC(e)$ 
according to Theorem \ref{srk-inf}, \cite[Theorem 3.7]{Pre} and \cite[Lemma 4.1]{Jan2}. 
We have therefore $(\xi\circ\iota)^{-1}(\mcO_{rmin}^{\SL_{n(2)}})$ is open in $\msN \cap \mfz(e)$.
\end{proof}

Since $\mcO_{e}\cap\mfz(e)$ is open in $\msN\cap \mfz(e)$, Lemma \ref{non-empty} ensures that
\begin{align*}
(\xi\circ\iota)^{-1}(\mcO_{rmin}^{\SL_{n(2)}})\bigcap \big(\mcO_{e}\cap\mfz(e)\big)\neq 
\emptyset. 
\end{align*}
We pick $ e_{0}\in \big(\mcO_{e}\cap \mfz(e)\big)
\bigcap(\xi\circ\iota)^{-1}(\mcO_{rmin}^{\SL_{n(2)}})$, then $(e,e_{0})=(\xi\circ \iota)(e_{0})\in 
\mcO_{rimn}^{\SL_{n(2)}}$. It follows that $\srk(\SL_{n(2)})=r_{(e,e_{0})}^{\SL_{n(2)}}$ 
by Theorem \ref{srk-inf}.

\begin{thm}{\label{height 2}}
Keep the notations for $\mcG, e$ and the assumption for $p$ in this section. Then we have
$\srk(\SL_{n(2)})=2(n-1)$.
\end{thm}
\begin{proof}
Since $\srk(\SL_{n(1)})=n-1$, we only need to prove $\srk(\SL_{n(2)})\geq 2(n-1)$ according to
Lemma \ref{inequa}.
Let $a=(a_{1},\ldots,a_{n-1})\in \mbG_{a}^{\times(n-1)}$ be a $(n-1)$-tuple and 
$e_{a}=\sum_{i=1}^{n-1}a_{i}e^{i}$. We define a closed subgroup of $\SL_{n}$, that is 
$U_{e}:=\langle 1+ e_{a}\mid a\in \mbG_{a}^{\times (n-1)} \rangle$ generated by $1+e_{a}$ when $a$ ranges over all
elements in $\mbG_{a}^{\times(n-1)}$. Additionally $U_{e}$ is an abelian unipotent 
subgroup of $\SL_{n}$. We denote by $\mfu_{e}:=\mrlie(U_{e})$ the abelian Lie algebra of $U_{e}$ and by
$\xymatrix{\phi: U_{e}\ar@{^{(}->}[r] & \SL_{n}}$ the closed embedding with 
the associated morphism $\xymatrix{d\phi: \mfu_{e}\ar@{->}[r] & \mfsl_{n}(\Bk)}$ of corresponding
restricted Lie algebras. One can easily check that for any $\Bk$-algebra $A$ and every $p$-nilpotent
element $x\in \mfu_{e}\otimes_{\Bk}A$ the homomorphism 
$\xymatrix{\exp_{d\phi(x)}: \mbG_{a}\otimes_{\Bk}A \ar@{->}[r] & \SL_{n}\otimes_{\Bk}A}$ factors through
$U_{e}\otimes_{\Bk}A$. Therefore, $\phi$ is an embedding of exponential type. By \cite[Lemma 1.7]{SFB1},
this gives $V_{2}(U_{e})=\set{(\alpha_{0},\alpha_{1})\in \mfu_{e}\times \mfu_{e}}=\mfu_{e}\times \mfu_{e}$. 
We observe that $(e,e_{0})\in V_{2}(U_{e})$. Then
$\xymatrix{\exp_{\underline{(e,e_{0})}}: \mbG_{a(2)}\ar@{->}[r] & \SL_{n(2)}}$ factors through $U_{e(2)}$ the 
second-Frobenius kernel of $U_{e}$.
Take the unique elementary abelian subgroup $\mcE_{U_{e(2)}}\subseteq U_{e(2)}$ given by 
\cite[Lemma 6.2.1]{Far4}, then the image of map $\exp_{\underline{(e,e_{0})}}$ will be contained in
$\mcE_{U_{e(2)}}$. Notice that $\cx_{\mcE_{U_{e(2)}}}(\Bk)=\cx_{U_{e(2)}}(\Bk)=\dim V_{2}(U_{e})=2 \dim \mfu_{e}=2(n-1)$, this gives
 $\srk(\SL_{n(2)})=r_{(e,e_{0})}^{\SL_{n(2)}}\geq 2(n-1)$, as desired.
\end{proof}

\newcommand{\etalchar}[1]{$^{#1}$}

\end{document}